\documentclass{article}
\usepackage{amsmath}
\usepackage{amsthm}
\usepackage{amsfonts}
\usepackage{amscd}

\newcommand{\shrinkmargins}[1]{
  \addtolength{\textheight}{#1\topmargin}
  \addtolength{\textheight}{#1\topmargin}
  \addtolength{\textwidth}{#1\oddsidemargin}
  \addtolength{\textwidth}{#1\evensidemargin}
  \addtolength{\topmargin}{-#1\topmargin}
  \addtolength{\oddsidemargin}{-#1\oddsidemargin}
  \addtolength{\evensidemargin}{-#1\evensidemargin}
  }
\shrinkmargins{.7}

\newcommand{\mrm}{\mathrm}
\newcommand{\mf}{\mathbf}

\newcommand{\ad}{{\mf{A}}}
\newcommand{\cc}{{\mf{C}}}
\newcommand{\zz}{{\mf{Z}}}
\newcommand{\qq}{{\mf{Q}}}
\newcommand{\rr}{{\mf{R}}}
\newcommand{\pp}{{\mf{P}}}
\newcommand{\ff}[1]{{\mf F}_{\!#1}}
\newcommand{\htt}{{\mf{T}}}

\newcommand{\ga}{\gamma}
\newcommand{\Ga}{\Gamma}
\newcommand{\eps}{\epsilon}

\newcommand{\rar}{\rightarrow}
\newcommand{\nin}{{\not\in}}
\newcommand{\into}{\hookrightarrow}
\newcommand{\dx}{{\,dx}}
\newcommand{\dy}{{\,dy}}
\newcommand{\uh}{{\mathfrak{H}}}
\newcommand{\cO}{{\mathcal{O}}}

\newcommand{\SL}{\mrm{SL}}
\newcommand{\PSL}{\mrm{PSL}}
\newcommand{\GL}{\mrm{GL}}

\newcommand{\jacobi}[2]{{\left(\frac{#1}{#2}\right)}}
\newcommand{\zmod}[1]{\zz/{#1}\zz}
\newcommand{\zmods}[1]{{(\zmod{#1})^*}}
\newcommand{\fdef}[3]{{{#1}\!:{#2}\rar{#3}}}
\newlength{\setoftmpheight}
\newcommand{\setof}[2]{{\settoheight{\setoftmpheight}{${#2}$}\left\{
  \left.{#1}{\vrule width 0pt height \setoftmpheight}\ \right|{#2}\right\}}}
\newcommand{\bigparens}[1]{{\left(\vbox{\vskip 2pt \hbox{${#1}$}}\right)}}
\newcommand{\smbigparens}[1]{{\left(\vbox{\vskip 1pt \hbox{${#1}$}}\right)}}
\newcommand{\bracket}[1]{{\left<{#1}\right>}}
\newcommand{\narpmod}[1]{{\left(\mrm{mod}\,{#1}\right)}}

\newcommand{\mattwo}[4]{{\begin{pmatrix}{#1}&{#2}\\{#3}&{#4}\end{pmatrix}}}
\newcommand{\smallmattwo}[4]{{\left(\begin{smallmatrix}{#1}&{#2}\\
  {#3}&{#4}\end{smallmatrix}\right)}}

\newcommand{\expn}[1]{{e^{2\pi\sqrt{-1}#1}}}

\newcommand{\ado}{{\ad^\infty}}
\newcommand{\zzhat}{{\hat{\zz}}}

\newcommand{\fff}{{\mf{f}}}
\newcommand{\ggg}{{\mf{g}}}

\newcommand{\Xis}[1]{{X_{\simeq,#1}}}
\newcommand{\Xise}{\Xis{\eps}}
\newcommand{\XisN}{{X_{\simeq}(N)}}
\newcommand{\XiseN}{{\Xise(N)}}

\newcommand{\Mt}{{\mathrm{M}_2}}
\newcommand{\SLt}[1]{{\SL_2(\zmod{#1})}}
\newcommand{\SLtp}{{\SLt{p}}}
\newcommand{\SLtN}{{\SLt{N}}}
\newcommand{\GLt}[1]{{\GL_2(\zmod{#1})}}

\newcommand{\GLtN}{{\GLt{N}}}
\newcommand{\GLtpQ}{{\GL_2^+(\qq)}}
\newcommand{\GLtpR}{{\GL_2^+(\rr)}}

\newcommand{\GLtA}{{\GL_2(\ado)}}
\newcommand{\GLtAA}{{\GLtA\times\GLtA}}
\newcommand{\GLtZh}{{\GL_2(\zzhat)}}
\newcommand{\GLtZZh}{{\GLtZh\times\GLtZh}}
\newcommand{\GLtQp}{{\GL_2(\qq_p)}}

\newcommand{\GLtZp}{{\GL_2(\zz_p)}}

\newcommand{\PSLt}[1]{{\PSL_2(\ff{#1})}}
\newcommand{\PSLtp}{{\PSLt{p}}}

\newcommand{\Goo}{{\Ga(1)\times\Ga(1)}}
\newcommand{\GwNN}{{\Ga_w(N)\times\Ga_w(N)}}

\newcommand{\Gis}[1]{{\Gamma_{\simeq,#1}}}
\newcommand{\Gise}{\Gis{\eps}}

\newcommand{\Diss}[1]{{\Delta_{\simeq,#1}^*}}
\newcommand{\Dises}{\Diss{\eps}}
\newcommand{\Disees}{\Diss{\eps,\eps'}}

\newcommand{\uhh}{{\uh\times\uh}}
\newcommand{\uhhst}{{\uh^*\times\uh^*}}
\newcommand{\uhG}{{\uh\times\GLtA}}
\newcommand{\uhhGG}{{\uhh\times\GLtAA}}

\newcommand{\htts}{\htt^*}      % No extra braces.
\newcommand{\httb}{{\overline{\htt}}}
\newcommand{\htteqs}{{\htts_\equiv}}
\newcommand{\httsarg}[2]{{\htts_{{#1},{#2}}}}
\newcommand{\httks}[1]{{\httsarg{k}{#1}}}
\newcommand{\httkes}{{\httks{\eps}}}
\newcommand{\httbarg}[2]{{\httb_{{#1},{#2}}}}
\newcommand{\httkb}[1]{{\httbarg{k}{#1}}}
\newcommand{\httkeb}{{\httkb{\eps}}}
\newcommand{\httiss}[1]{{\htts_{{#1},\simeq}}}
\newcommand{\httkiss}{{\httiss{k}}}
\newcommand{\httisb}[1]{{\httb_{{#1},\simeq}}}
\newcommand{\httkisb}{{\httisb{k}}}

\newcommand{\PoQ}{{\pp^1(\qq)}}

\newcommand{\pairmattwo}[8]{\bigparens{
  \mattwo{#1}{#2}{#3}{#4},\mattwo{#5}{#6}{#7}{#8}}}
\newcommand{\pairsmallmattwo}[8]{\smbigparens{
  \smallmattwo{#1}{#2}{#3}{#4},\smallmattwo{#5}{#6}{#7}{#8}}}
\newcommand{\genpairmattwo}{\pairmattwo{a_1}{b_1}{c_1}{d_1}
  {a_2}{b_2}{c_2}{d_2}}

\newcommand{\qmp}{{\qq(\sqrt{-p})}}

\newcommand{\Sbar}{{\overline{S}}}
\newcommand{\Sk}[2]{{S_k(\Gis{#1}({#2}))}}
\newcommand{\Ske}[1]{{\Sk{\eps}{#1}}}
\newcommand{\SkN}[1]{{\Sk{#1}{N}}}
\newcommand{\SkeN}{\SkN{\eps}}
\newcommand{\Skb}[2]{{\Sbar_k(\Gis{#1}({#2}))}}
\newcommand{\Skeb}[1]{{\Skb{\eps}{#1}}}
\newcommand{\SkNb}[1]{{\Skb{#1}{N}}}
\newcommand{\SkeNb}{\SkNb{\eps}}
\newcommand{\Skis}[1]{{S_{k,\simeq}({#1})}}
\newcommand{\SkisN}{{\Skis{N}}}
\newcommand{\Skisb}[1]{{\Sbar_{k,\simeq}({#1})}}
\newcommand{\SkisNb}{{\Skisb{N}}}
\newcommand{\Stgis}[2]{{S_{(2,2)}(\Gis{#1}({#2}))}}
\newcommand{\Stgisp}[1]{{\Stgis{#1}p}}
\newcommand{\Step}{\Stgisp{\eps}}
\newcommand{\Stepb}{{\Sbar_{(2,2)}(\Gise(p))}}
\newcommand{\Stmp}{{\Stgis{-1}{p}}}
\newcommand{\Sistp}{{S_{(2,2),\simeq}(p)}}
\newcommand{\Sistpb}{{\Sbar_{(2,2),\simeq}(p)}}
\newcommand{\SGw}[2]{{S_{#1}(\Ga_w({#2}))}}
\newcommand{\SGwN}[1]{{\SGw{#1}{N}}}
\newcommand{\SkGwN}{{\SGwN{k}}}
\newcommand{\SkkGwN}{{\SGwN{k_1}\otimes \SGwN{k_2}}}
\newcommand{\SGwb}[2]{{\Sbar_{#1}(\Ga_w({#2}))}}
\newcommand{\SGwNb}[1]{{\SGwb{#1}{N}}}
\newcommand{\StGwp}{{\SGw{2}{p}}}
\newcommand{\SkUN}{{S_k(U(N))}}
\newcommand{\StUp}{{S_2(U(p))}}
\newcommand{\SkkGwNb}{{\SGwNb{k_1}\otimes \SGwNb{k_2}}}

\newcommand{\KNb}[2]{{\overline{K}_{{#1},{#2}}(N)}}
\newcommand{\KeNb}[1]{{\KNb{#1}{\eps}}}
\newcommand{\KkNb}[1]{{\KNb{k}{#1}}}
\newcommand{\KkeNb}{{\KkNb{\eps}}}
\newcommand{\KkNbp}[1]{{\overline{K}'_{k,{#1}}(N)}}
\newcommand{\KkeNbp}{{\KkNbp{\eps}}}
\newcommand{\Kis}{{K_\simeq}}

\newcommand{\diam}[1]{{\langle{#1}\rangle}}

\newcommand{\Tmi}{{T_{m_1,m_2}}}
\newcommand{\Tni}{{T_{n_1,n_2}}}

\newcommand{\cmi}{{c_{m_1,m_2}}}
\newcommand{\cni}{{c_{n_1,n_2}}}
\newcommand{\dmi}{{d_{m_1,m_2}}}

\newcommand{\lmi}{{\lambda_{m_1,m_2}}}
\newcommand{\lni}{{\lambda_{n_1,n_2}}}

\newcommand{\Uis}[1]{{U_{\simeq}({#1})}}
\newcommand{\UisN}{{\Uis{N}}}

\newcommand{\Sigb}{\overline{\Sigma}}

\newtheorem{theorem}{Theorem}[section]
\newtheorem{prop}[theorem]{Proposition}
\newtheorem{lemma}[theorem]{Lemma}
\newtheorem{cor}[theorem]{Corollary}
\newtheorem{conj}[theorem]{Conjecture}

\numberwithin{equation}{section}

\begin{document}
\bibliographystyle{abbrv}
\title{Moduli for Pairs of Elliptic Curves with Isomorphic
  $N$-torsion}
\author{David Carlton}

\maketitle

\begin{abstract}
  We study the moduli surface for pairs of elliptic curves together
  with an isomorphism between their $N$-torsion groups.  The Weil
  pairing gives a ``determinant'' map from this moduli surface to
  $\zmods N$; its fibers are the components of the surface.  We define
  spaces of modular forms on these components and Hecke
  correspondences between them, and study how those spaces of modular
  forms behave as modules for the Hecke algebra.  We discover that the
  component with determinant $-1$ is somehow the ``dominant'' one; we
  characterize the difference between its spaces of modular forms and
  the spaces of modular forms on the other components using forms with
  complex multiplication.  Finally, we show some simplifications that
  arise when $N$ is prime, including a complete determination of such
  CM-forms, and give numerical examples.
\end{abstract}

\tableofcontents

\section{Introduction}

If $R$ is the ring of integers in a totally real number field, one can
consider the Hilbert modular variety associated to $R$, which
parameterizes abelian varieties of dimension $[R: \zz]$ together with
a map from $R$ into their endomorphism ring.  This modular variety is
disconnected; its components correspond to polarization types, and are
indexed by elements of the narrow class group of $R$.  One can define
spaces of modular forms associated to the modular variety and to its
components; the former are more adelic in nature, while the latter are 
more classical.

In this paper, we consider a variant of the above situation, where we
replace $R$ by the order $(\zz \times \zz)_{\equiv (N)}$ that consists
of pairs of integers that are congruent mod $N$.  Thus, we replace
our totally real number field by the totally real ``number algebra''
$\qq \times \qq$, and in addition consider a non-maximal order rather
than the full ring of integers.  As in the traditional situation, one
can associate a modular variety to this situation, and study its
components, which are indexed by $\zmods{N}$; this has been done in
Hermann~\cite{hermann} and Kani and Schanz~\cite{kani-modular}.  One
can also define spaces of classical and adelic modular forms, which we
do in this paper.

These ``degenerate'' Hilbert modular varieties and modular forms
should have properties very similar to those of traditional Hilbert
modular varieties and modular forms.  However, they can also be
related to modular curves and elliptic modular forms, which have been
the subject of extensive study.  For example, these surfaces have an
interpretation as moduli spaces for pairs of elliptic curves with
isomorphic $N$-torsion, and can be constructed as a quotient of
$X(N)\times X(N)$.  Thus, we expect them to be a particularly suitable
test ground for exploring properties of Hilbert modular surfaces and
modular forms.  We expect the generalization to the case where $R$ is
an order in a product of ring of integers of totally real number
fields to be of interest as well: for example, $R$ might be the Hecke
algebra $\htt_0(N)$ associated to the modular curve $X_0(N)$.

One such new property, which is the main goal of this paper, involves
studying how these components of the degenerate Hilbert modular
variety vary.  It is easy to see that two components whose index
differs by a square are isomorphic, but there is no reason why other
components should be isomorphic.  Indeed, Hermann has shown that, for
example, if $N=7$ then the component indexed by 1 is a rational
surface and the component indexed by $-1$ is a K3 surface; similarly,
if $N = 11$, the component indexed by 1 is an elliptic surface and the
component indexed by $-1$ is of general type.  As Kani and Schanz
noted, this change in geometric complexity is reflected by the
geometric genera of the components.

We show in Section~\ref{sec:cusp-relationships} that, for $N$ fixed
and prime, the component indexed by $-1$ always has the largest
geometric genus of any of the components; we give an explicit formula
for the difference of geometric genera in Section~\ref{sec:prime}.
The geometric genus of a component is the dimension of a suitable
space of cusp forms; we exhibit this difference in genera as the
dimension of a certain special subspace of the space of cusp forms on
the $-1$ surface; we call it the \emph{Hecke kernel} since it can
be seen as the intersection of the kernels of certain Hecke operators.
We also show in Section~\ref{sec:hecke-kernel} that the elements of
the Hecke kernel have an alternative characterization as forms with
complex multiplication; we give an explicit construction of the forms
in Section~\ref{sec:examples}.  The proof of these results involves
the interplay between spaces of adelic and classical modular forms.

I would like to thank Fred Diamond, Jordan Ellenberg, Steven Kleiman,
and Barry Mazur for the help that they have given me while writing
this; and the N.D.S.E.G. Fellowship Program for the support that it
has provided.

\section{Basic Definitions}
\label{sec:basic}

Let $X_w(N)$ be the curve over $\cc$ parameterizing elliptic curves
together with a basis for their $N$-torsion that maps to some
specified $N$'th root of unity under the Weil pairing.\footnote{This
  curve is traditionally denoted by $X(N)$; however, we have chosen to
  use the notation $X(N)$ to denote the (geometrically reducible)
  curve coming from the adelic mod $N$ principal congruence subgroup,
  and have changed all notation accordingly.}  It is Galois over the
curve $X_w(1)$ with Galois group $\SLtN/\{\pm 1\}$.  Let $\SLtN$ act
on the product surface $X_w(N)\times X_w(N)$ via the diagonal action;
we can then form the quotient surface, which we shall denote by
$\Xis{1}(N)$.  More generally, if $\eps$ is an element of $\zmods N$
and if $\SLtN$ acts on the first factor via the natural action but on
the second factor via the automorphism
\[
\theta_\eps\!: \mattwo abcd \mapsto \mattwo a{\eps^{-1} b}{\eps c}d
\]
then we denote the quotient surface by $\XiseN$.  And we set
\[
\XisN = \coprod_{\eps \in \zmods{N}} \Xise(N).
\]

These surfaces can also be constructed in another fashion, as
degenerate Hilbert modular surfaces: let $\uh$ be the upper half
plane, with $\Ga(1) = \SL_2(\zz)$ acting on it via fractional linear
transformations.  Then $\Goo$ acts on $\uhh$; if we denote by
$\Gise(N)$ the subgroup of $\Goo$ given by
\[
\setof{\genpairmattwo}{
\begin{array}{rcll}
  a_1 & \equiv & a_2 & \pmod N,\\
  b_1 & \equiv & \eps b_2 & \pmod N, \\
  \eps c_1 & \equiv & c_2 & \pmod N, \\
  d_1 & \equiv & d_2 & \pmod N
\end{array}
}
\]
then the quotient $\Gise(N) \backslash \uhh$ is an open subset of
$\XiseN$, and if we denote by $\uh^*$ the space $\uh \coprod \PoQ$
then $\Gise(N) \backslash \uhhst$ is all of $\XiseN$.

The surface $\XiseN$ (or, more properly, the open subset given by
using $\uhh$ instead of $\uhhst$) is a coarse moduli space for triples
$(E_1, E_2, \phi)$ where the $E_i$'s are elliptic curves and $\phi$ is
an isomorphism from $E_1[N]$ to $E_2[N]$ such that $\wedge^2 \phi$
raises the Weil pairing to the $\eps$'th power.  The modular
parameterization is given as follows: let $(\tau_1,\tau_2) \in \uhh$
and let $E_i$ be the elliptic curve given by the lattice with basis
$\{1,\tau_i\}$.  Also, let $e$ be an integer that reduces to $\eps$
mod $N$.  We then have the map $\phi$ from $E_1[N]$ to $E_2[N]$ that
sends $\tau_1/N$ to $e\tau_2/N$ and $1/N$ to $1/N$; it raises the Weil
pairing to the $\eps$'th power, the group of elements of $\Goo$ that
preserve $\phi$ is the subgroup $\Gise(N)$ defined above, and every
triple $(E_1, E_2, \phi)$ arises in this fashion.

The structure of the $\XiseN$'s as complex surfaces has been studied
by Hermann in \cite{hermann} and by Kani and Schanz in
\cite{kani-modular}; our $\XiseN$ is Hermann's $Y_{N,\eps^{-1}}$ and
Kani and Schanz's $Z_{N,\eps^{-1}}$.\footnote{We replaced their $\eps$
  by $\eps^{-1}$ to simplify the normalizations given in
  Section~\ref{sec:adelic}; since $\XiseN$ and $\Xis{\eps^{-1}}(N)$
  are isomorphic, this is an unimportant change.}  In particular, Kani
and Schanz give explicit formulas and tables computing various
invariants of the $\XiseN$'s, such as the dimensions of various
cohomology groups.  They also give explicit minimal desingularizations
of the surfaces.

We now define spaces of modular forms on these surfaces.  Thus, let
$\fdef f\uhh\cc$ be a holomorphic function; let $\ga = (\ga_1,\ga_2)$
be an element of $\GLtpR\times\GLtpR$, where $\GLtpR$ is the set of
elements of $\GL_2(\rr)$ with positive determinant; and let $k =
(k_1,k_2)$ be a pair of natural numbers.  We define the function
$\fdef{f|_{k,\ga}}\uhh\cc$ by
\[
f|_{k,\ga}(z_1,z_2)= f(\ga_1(z_1),\ga_2(z_2))j(\ga_1,z_1)^{-k_1}
  j(\ga_2,z_2)^{-k_2}
\]
where, if $\sigma = \smallmattwo abcd$ is an element of $\GLtpR$, then
$\sigma(z) = (az + b)/(cz+d)$ and
\[
j(\sigma,z) = (ad-bc)^{-1/2}(cz+d).
\]
We write $f|_\ga$ instead of $f|_{k,\gamma}$ if $k$ is clear from
context.

Defining $\Ga(1)$ to be $\SL_2(\zz)$, we say that a subgroup $\Ga$ of
$\Goo$ is a \emph{congruence subgroup} if it contains the group
$\GwNN$ for some $N$, where $\Ga_w(N)$ is defined to be the set of
matrices in $\SL_2(\zz)$ that are congruent to the identity mod $N$.
A function $\fdef f\uhh\cc$ is a \emph{modular form for $\Ga$ of
  weight $k$} if $f|_{k,\ga} = f$ for all $\ga \in \Ga$ and if $f$ is
holomorphic at the cusps.  To explain this latter condition, assume
that $\GwNN \subset\Ga$.  Then $f(z_1+N,z_2) = f(z_1,z_2)$ for all
$(z_1,z_2) \in \uhh$; so setting $q_1 = \expn{z_1/N}$, we can write
\[
f(z_1,z_2) = \sum_{m \in\zz} c_m(f)(z_2)q_1^m
\]
for some functions $c_m(f)$.  If $c_m(f)$ is zero for all $m < 0$ and
if a similar condition holds if we do a Fourier expansion in $z_2$, we
say that $f$ is \emph{holomorphic at infinity}.  And $f$ is
\emph{holomorphic at all of the cusps} if, for all $\ga \in \Goo$,
$f|_{k,\ga}$ is holomorphic at infinity.

A modular form is a \emph{cusp form} if it vanishes at all of the
cusps; that is to say, if whenever we take a Fourier expansion of
$f|_{k,\ga}$ in either variable as above, $c_0(f)$ is zero.  We denote
the space of all modular forms of weight $k$ for $\Ga$ by $M_k(\Ga)$;
we denote the space of all cusp forms by $S_k(\Ga)$.

If $\Ga = \Ga_1 \times \Ga_2$, with each $\Ga_i$ a congruence subgroup
of $\Ga(1)$, then there is a natural map from $M_{k_1}(\Ga_1) \otimes
M_{k_2}(\Ga_2)$ to $M_{(k_1,k_2)}(\Ga_1\times\Ga_2)$ which sends $f_1
\otimes f_2$ to the function
\[
(z_1,z_2) \mapsto f_1(z_1)f_2(z_2).
\]
Furthermore, this map sends cusp forms to cusp forms.  It is in fact
an isomorphism in either the modular form or cusp form case:

% Local notation, only to be used in this proposition and its proof.
\newcommand{\MkkGG}[1]{{M_{(k_1,k_2)}(\Ga_1\times\Ga_2,#1)}}
\newcommand{\MkkGGSS}{\MkkGG{(S_1\times\uh^*)\cup(\uh^*\times S_2)}}

\begin{prop}
  \label{prop:group-prod}
  If $S$ is a subset of $\uh^*$ or $\uhhst$ and $\Ga$ is a congruence
  subgroup of $\Ga(1)$ or $\Gamma(1)\times\Gamma(1)$, let $M_k(\Ga,S)$
  be the set of forms in $M_k(\Ga)$ that vanish on the points in $S$.
  Then for any congruence subgroups $\Ga_1$ and $\Ga_2$ of $\Ga(1)$
  and subsets $S_1$ and $S_2$ of $\uh^*$, the natural map
  \[
  M_{k_1}(\Ga_1,S_1) \otimes M_{k_2}(\Ga_2,S_2) \rar \MkkGGSS
  \]
  is an isomorphism.
\end{prop}

\begin{proof}
  This follows by induction on the dimension of
  $M_{k_1}(\Gamma_1,S_1)$.
\end{proof}

\begin{cor}
  \label{cor:group-prod}
  Given any natural numbers $k_1$, $k_2$, and $N$, we have
  isomorphisms
  \[
  M_{(k_1,k_2)}(\Gise(N)) = (M_{k_1}(\Ga_w(N))\otimes
  M_{k_2}(\Ga_w(N)))^{\SLtN}
  \]
  and
  \[
  S_{(k_1,k_2)}(\Gise(N)) = (S_{k_1}(\Ga_w(N))\otimes
  S_{k_2}(\Ga_w(N)))^{\SLtN},
  \]
  where $\SLtN$ acts on the first member of the tensor product in the
  natural fashion and on the second member via the automorphism
  $\theta_\eps$.
\end{cor}

\begin{proof}
  By Proposition~\ref{prop:group-prod},
  \[
  M_{(k_1,k_2)}(\Ga_w(N)\times\Ga_w(N)) = (M_{k_1}(\Ga_w(N))\otimes
  M_{k_2}(\Ga_w(N)));
  \]
  that $\SLtN$-invariants correspond to forms in
  $M_{(k_1,k_2)}(\Gise(N))$ follows from the definitions.  The cusp
  form case is similar, setting $S_1$ and $S_2$ in the Proposition to
  be equal to $\pp^1(\qq)$.
\end{proof}

This allows us to express the dimension of the space
$S_{(2,2)}(\Gise(N))$ in terms of data given in Kani and
Schanz~\cite{kani-modular}:

\begin{cor}
  \label{cor:cusp-dim-coh}
  The dimensions of the spaces $S_{(2,2)}(\Gise(N))$ and
  $H^2(\XiseN,\cO_\XiseN)$ are equal, and they are also equal to the
  geometric genus of a desingularization of $\XiseN$.
\end{cor}

\begin{proof}
  We have the equalities
  \begin{align*}
    \dim S_{(2,2)}(\Gise(N)) &= \dim (S_2(X_w(N)) \otimes
    S_2(X_w(N)))^{\SLtN} \\
    &= \dim (H^1(X_w(N),\cO_{X_w(N)})\otimes
    H^1(X_w(N),\cO_{X_w(N)}))^{\SLtN} \\
    &= \dim H^2(X_w(N)\times X_w(N),\cO_{X\times X})^{\SLtN} \\
    &= \dim H^2(\SLtN \backslash (X_w(N)\times X_w(N)),
    \cO_{\SL_2\backslash X\times X}) \\
    &= \dim H^2(\XiseN,\cO_\XiseN).
  \end{align*}
  This last quantity is equal to the geometric genus, by Kani and
  Schanz~\cite{kani-surfaces}, Proposition~3.1.
\end{proof}

Of course, this isn't too surprising: weight 2 cusp forms should
correspond to holomorphic 2-forms.

If $f$ is a modular form on $\Gise(N)$, it
has a Fourier expansion
\[
f(z_1,z_2) = \sum_{m_1,m_2 \ge 0} \cmi(f) q_1^{m_1}q_2^{m_2}
\]
where $q_i = \expn{z_i/N}$.  There is one thing that we can say
immediately about the Fourier coefficients $\cmi(f)$:

\begin{prop}
  \label{prop:zero-coeff}
  For all $f \in M_{(k_1,k_2)}(\Gise(N))$, the Fourier coefficient
  $\cmi(f)$ is zero unless $\eps m_1 + m_2 \equiv 0 \pmod N$.
\end{prop}

\begin{proof}
  This follows from the fact that $f = f|_{\pairsmallmattwo
    1e011101}$, where $e$ is an integer congruent to $\eps$ mod $N$.
\end{proof}

Thus, most of the Fourier coefficients are ``missing''.  This turns
out to make it natural to also study modular forms on the surface
$\XisN$, even when we are only interested in one of the individual
$\XiseN$'s; we shall elaborate on this theme in
Section~\ref{sec:cusp-relationships}.

One way to produce forms on $\XiseN$ is to consider forms on
$\Xise(N/d)$ to be forms on $\Xise(N)$, for $d$ a divisor of $N$.
Such forms have Fourier coefficients $c_{m_1,m_2}$ equal to zero
unless $d$ divides $m_1$ (and hence $m_2$, by
Proposition~\ref{prop:zero-coeff}).  The converse is also true:

\begin{theorem}
  \label{thm:level-change}
  Let $f$ be a modular form of weight $k$ on $\Gise(N)$,
  and assume that, for some $d | N$, we have $c_{m_1,m_2}(f) = 0$
  unless $d|m_1$.  Then $f$ is an element of $M_k(\Gise(N/d))$.
\end{theorem}

\begin{proof}
  The fact that $c_{m_1,m_2}(f) = 0$ unless $d|m_1$ is equivalent to
  having $f$ be invariant under
  \[
  \pairmattwo1{N/d}011001.
  \]
  Thus, we have to show that the smallest subgroup $\Ga$ containing
  both $\pairsmallmattwo1{N/d}011001$ and $\Gise(N)$ is $\Gise(N/d)$.
  Furthermore, we can take the quotient by $\Ga_w(N) \times \Ga_w(N)$,
  and thus consider all matrices to be elements of $\SLtN$.  Letting
  $G = \setof{\gamma \in \SLtN}{(\gamma,1) \in \Gamma}$, we see that
  $\Gamma = G\times\{1\}\cdot \Gise(N)$ and that $\Gamma$ is a subgroup
  if and only if $G$ is normal.
  Thus, we have to show that the smallest normal subgroup of $\SLtN$
  containing the matrix $\tau_{N/d} = \smallmattwo 1{N/d}01$ is the
  kernel of the natural map from $\SLtN$ to $\SLt{(N/d)}$.
  Furthermore, we can assume that $d$ is a prime $p$, and by the
  Chinese remainder theorem we can assume that $N=p^l$ for some $l$.
  
  First, assume that $l = 1$, so we want to show that the smallest
  normal subgroup $G$ of $\SLtp$ containing $\tau_1 = \smallmattwo
  1101$ is the entire group.  We first look at the image of $G$ in
  $\PSL_2(\zmod{p})$.  If $p > 3$ then $\PSL_2(\zmod{p})$ is simple,
  so the image of $G$ is all of $\PSL_2(\zmod{p})$.  If $p = 3$ then
  $\PSL_2(\zmod{3})$ is isomorphic to $A_4$ and $\tau_1$ is an element
  of order 3; but since the only proper normal subgroups of $A_4$
  contain only elements of order 1 and 2, we again have that the image
  of $G$ is all of $\PSL_2(\zmod{3})$.  Similarly, if $p=2$, then
  $\PSL_2(\zmod{2})$ is isomorphic to $S_3$ and $\tau_1$ has order 2,
  so again our image must be all of $\PSL_2(\zmod 2)$.
  
  This implies that $G$ must either be all of $\SLtp$ or a subgroup of
  index two which projects onto all of $\PSL_2(\zmod p)$.  But if $p =
  2$ then $\SLt{2} = \PSL_2(\zmod 2)$; if $p = 3$ then $\SLt{3}$ has
  only two non-trivial one-dimensional representations, whose kernels
  are of index 3; and if $p > 3$ then $\SLtp$ has no non-trivial
  one-dimensional representations, so again has no subgroups of index
  2.
  
  Finally, assume that $l > 1$, and that we have a normal subgroup $G$
  containing $\tau_q$, where $q=p^{l-1}$.  (Note that $q^2$ is zero in
  $\zmod{p^l}$, which greatly simplifies calculations.)  We then have
  to show that $G$ contains all matrices of the form
  $\smallmattwo{1+aq}{bq}{cq}{1+dq}$ with determinant 1; this
  condition on the determinant is equivalent to having $a$ equal to
  $-d$ in $\zmod p$.  But it is easy to produce all such matrices by
  taking suitable multiples of $\tau_q$, its conjugate by
  $\smallmattwo 0{-1}10$, and its conjugate by $\smallmattwo a{a-1}11$.
\end{proof}

We hope that the following stronger result is true:

\begin{conj}
  \label{conj:level-change}
  Let $f$ be a modular form on $\Gise(N)$ such that $c_{m_1,m_2}(f) =
  0$ unless $(m_i,N) > 1$.  Then $f$ can be written as a sum of
  modular forms $f_j$ on $\Gise(N/p_j)$ where the $p_j$'s are the
  prime divisors of $N$.  Furthermore, if $f$ is a cusp form then the
  $f_j$ can be chosen to be cusp forms.
\end{conj}

Of course, Theorem~\ref{thm:level-change} implies
Conjecture~\ref{conj:level-change} for $N$ a prime power.  They are
both analogous to results proved as parts of Atkin-Lehner theory on
the curves $X_1(N)$.  (C.f. Theorem~1 of
Atkin-Lehner~\cite{atkin-lehner} or Lang~\cite{lang-mod-forms},
Theorem~VIII.3.1.)

We let $\SkeNb$ be the quotient of $S_k(\Gise(N))$ by the subgroup of
forms $f$ whose Fourier coefficients $c_{m_1,m_2}(f)$ are zero unless
$(m_i,N) > 1$.  In the $X_1(N)$ case, this would have the effect of
replacing $S_k(\Ga_1(N))$ by a space with the same Hecke eigenspaces
but where each eigenspace is one-dimensional, generated by the newform
in that eigenspace; we shall see in Theorem~\ref{thm:cusp-free} that
Hecke eigenspaces in $\SkeNb$ are also one-dimensional.  Finally, we
let
\[
\SkisN = \prod_{\eps \in \zmods{N}} S_k(\Gise(N)),
\]
and we let
\[
\SkisNb = \prod_{\eps\in\zmods{N}} \SkeNb.
\]
Note that in the definitions of $\SkeNb$ and $\SkisNb$ it's enough to
assume that the Fourier coefficients are zero unless $(m_1,N) > 1$ (or
unless $(m_2,N) > 1$), by Proposition~\ref{prop:zero-coeff}.

\begin{prop}
  \label{prop:p-nobar}
  The spaces $\Step$ and $\Stepb$ are equal, as are the spaces
  $\Sistp$ and $\Sistpb$.
\end{prop}

\begin{proof}
  We have to show that if $f$ is an element of $\Step$ such that
  $c_{m_1,m_2}(f) = 0$ unless $p | m_1$ then $f$ is zero.
  Theorem~\ref{thm:level-change} implies that such an $f$ is in fact a
  form on $\Gise(1)$.  By Corollary~\ref{cor:group-prod}, $f$ can be
  considered to be an element of $S_2(\Ga(1)) \otimes S_2(\Ga(1))$.
  But $S_2(\Ga(1))$ is zero, so $f$ is zero.
\end{proof}

\begin{prop}
  \label{prop:bar-nobar}
  If $p$ is a prime then
  \[
  \dim \Ske{p^l} = \sum_{j=0}^l \dim \Skeb{p^j}.
  \]
\end{prop}

\begin{proof}
  This follows immediately from Theorem~\ref{thm:level-change}.
\end{proof}

Conjecture~\ref{conj:level-change} would imply a similar statement for
forms of arbitrary level.

\section{Hecke Operators on $\XiseN$}
\label{sec:hecke-one}

Set
\[
\Dises(N) = \setof{\genpairmattwo}{
  \begin{array}{l}
    a_i,b_i,c_i,d_i \in \zz, \\
    a_id_i-b_ic_i > 0, \\
    (a_id_i-b_ic_i,N) = 1, \\
    \begin{array}{rcll}
      a_1 & \equiv & a_2 & \pmod N,\\
      b_1 & \equiv & \eps b_2 & \pmod N, \\
      \eps c_1 & \equiv & c_2 & \pmod N, \\
      d_1 & \equiv & d_2 & \pmod N
    \end{array}
  \end{array}
  }.
\]
We can partition $\Dises(N)$ into double $\Gise(N)$-cosets; each
double coset is called a \emph{Hecke operator}.  They act on the
spaces of modular forms as follows:

Let $\ga = (\ga_1,\ga_2)$ be an element of $\Dises(N)$, and let
\[
\Gise(N)\ga\Gise(N) = \coprod_j \Gise(N)\ga_j
\]
be a decomposition of the double coset generated by $\ga$ into left
cosets.  Then for a form $f$ in $M_{(k_1,k_2)}(\Gise(N))$, we define
\[
f|_{(k_1,k_2),\Gise(N)\ga\Gise(N)} =
  \det(\ga_1)^{(k_1/2)-1}\det(\ga_2)^{(k_2/2)-1}\sum_j
  f|_{(k_1,k_2),\ga_j}.
\]
We see as in Shimura~\cite{shimura}, Chapter~3, that
$f|_{(k_1,k_2),\Gise(N)\ga\Gise(N)}$ is an element of the space
$M_{(k_1,k_2)}(\Gise(N))$, that cusp forms are transformed into cusp
forms, and that the product of two Hecke operators is a sum of Hecke
operators.

Let $\Tni$ be the operator given by the sum of the double
cosets containing elements $(\ga_1,\ga_2)$ where $\det(\ga_i) = n_i$.
This is zero unless $n_1 \equiv n_2 \pmod N$ and $(n_i,N) = 1$.  Left
coset representatives for it are given as follows:

\begin{prop}
  \label{prop:hecke-coset}
  Let $(n_1,n_2)$ be a pair of positive integers that are congruent
  mod $N$ and that are relatively prime to $N$.  The set of elements
  of $\Dises(N)$ that have determinant $(n_1,n_2)$ then has the
  following left coset decomposition:
  \[
  \coprod_{\substack{a_1,a_2>0\\a_id_i=n_i\\0\le b_i<d_i}}
  \Gise(N)\bigparens{\sigma_{a_1}\mattwo{a_1}{b_1N}{0}{d_1},
    \sigma_{a_2}\mattwo{a_2}{b_2N}{0}{d_2}}
  \]
  where, for $a \in \zmods{N}$, $\sigma_a$ is any matrix in $\Ga(1)$
  that is congruent to $\smallmattwo{a^{-1}}00a$ mod $N$.
\end{prop}

\begin{proof}
  First, note that the above cosets do indeed occur in $\Tni$.  Also,
  it is easy to see that the above cosets are disjoint.  Thus, we have 
  to show that the cosets cover all of $\Tni$.
  
  Let $(\delta_1,\delta_2)$ be an element of $\Dises(N)$ with
  determinant $(n_1,n_2)$.  By Shimura~\cite{shimura},
  Proposition~3.36, we can multiply $\delta_1$ on the left by an
  element of $\Ga(1)$ to get it into the form $\smallmattwo
  {a_1}{b_1}0{d_1}$, with $a_1>0$, $a_1d_1=n_1$, and $0 \le b_1 <
  d_1$.  Subsequently multiplying it on the left by an element of the
  form $\smallmattwo 1x01$ will put it into the form $\smallmattwo
  {a_1}{b_1N}0{d_1}$, but possibly with a different $b_1$.  (We can
  still force $b_1$ to be in the range $0 \le b_1 < d_1$, however.)
  And since $\sigma_{a_1}$ is an element of $\Ga(1)$, we have shown
  that there is an element $\ga_1$ of $\Ga(1)$ such that
  $\ga_1\delta_1$ is of the form
  $\sigma_{a_1}\smallmattwo{a_1}{b_1N}0{d_1}$.
  
  We can choose an element $\ga_2$ of $\Ga(1)$ such that
  $(\ga_1,\ga_2)$ is in $\Gise(N)$: reduce $\ga_1$ mod $N$, apply
  $\theta_\eps$ to it, and lift it back to $\Ga(1)$.  Multiplying
  $(\delta_1,\delta_2)$ on the left by $(\ga_1,\ga_2)$, we can thus
  assume that $\delta_1$ is of the form
  $\sigma_{a_1}\smallmattwo{a_1}{b_1N}0{d_1}$.  But then the
  congruence relations force $\delta_2$ to be congruent to the matrix
  \[
  \mattwo100{n_1} \equiv \mattwo100{n_2} \pmod N.
  \]
  
  Now that we have fixed $\delta_1$ to be of the correct form, we
  still have to force $\delta_2$ to be of the correct form, and we are
  only allowed to multiply $\delta_2$ on the left by elements of
  $\Ga_w(N)$.  Thus, we need to find an element $\ga'_2$ of $\Ga_w(N)$
  such that $\ga'_2\delta_2$ is of the form
  $\sigma_{a_2}\smallmattwo{a_2}{b_2N}0{d_2}$.  However, $\delta_2$ is
  in what Shimura calls $\Delta'$ (see Shimura~\cite{shimura}, p.~68),
  so we can indeed find such a $\ga'_2$ by Proposition~3.36 of
  Shimura~\cite{shimura}.
\end{proof}

The action of the Hecke operators $\Tni$ descends to the spaces
$\SkeNb$:

\begin{prop}
  \label{prop:hecke-equiv}
  If $f$ is a form in $S_k(\Gise(N))$ such that $c_{m_1,m_2}(f) = 0$
  unless $(N, m_i) > 1$ then $\Tni f$ has the same property for all
  $n_1 \equiv n_2 \pmod N$.
\end{prop}

\begin{proof}
  For $d | N$, define the operator $i_d$ by
  \[
  i_d(f) = \sum_{\substack{m_1,m_2 > 0\\d|m_1,m_2}}
  c_{m_1,m_2}(f)q_1^{m_1}q_2^{m_2};
  \]
  it has an equivalent definition as
  \[
  i_d(f) = \frac1d\sum_{0 \le e < d} f|_\pairsmallmattwo10011{Ne/d}01.
  \]
  By the principle of inclusion
  and exclusion, the statement that $c_{m_1,m_2}(f) = 0$ unless
  $(N,m_i) > 1$ is equivalent to having
  \[
  f = \sum_{p|N} i_p(f) - \sum_{\substack{p_1,p_2|N\\p_1<p_2}}
  i_{p_1p_2}(f) + \dotsb,
  \]
  and we want to show that if that is the case for $f$ then it is also
  the case for $\Tni f$.  It is therefore enough to show that
  $\Tni$ commutes with any $i_d$.  But
  \[
  \mattwo1{eN/d}01\sigma_{a_2}\mattwo{a_2}{b_2N}0{d_2}
  \]
  is congruent to
  \[
  \sigma_{a_2}\mattwo{a_2}{b_2N}0{d_2}\mattwo1{en_2N/d}01
  \]
  mod $N$, so by Proposition~\ref{prop:hecke-coset}, commuting with
  $\Tni$ simply permutes the $e$'s that occur in our alternate
  definition of $i_d$.
\end{proof}

Proposition~\ref{prop:hecke-equiv} would be an easy corollary to
Conjecture~\ref{conj:level-change}.

\begin{prop}
  \label{prop:self-adjoint}
  For all $(\delta_1,\delta_2) \in \Dises(N)$, the $\Gise(N)$-double
  cosets generated by $(\delta_1,\delta_2)$ and
  $(\delta_1^\iota,\delta_2^\iota)$ are equal, where
  \[
  {\mattwo abcd}^\iota = \mattwo d{-b}{-c}a.
  \]
\end{prop}

\begin{proof}
  We need to find matrices $(\ga_1,\ga_2)$ and $(\ga_1',\ga_2')$ in
  $\Gise(N)$ such that
  \[
  (\ga_1\delta_1,\ga_2\delta_2) =
  (\delta_1^\iota\ga_1',\delta_2^\iota\ga_2').
  \]
  Since $\delta_1$ and $\delta_1^\iota$ have the same elementary
  divisors, we can choose a $\ga_1$ and $\ga_1'$ that give us equality
  on the first coordinate.  Now pick $\ga_2$ and $\ga_2'$ such that
  $(\ga_1,\ga_2)$ and $(\ga_1',\ga_2')$ are in $\Gise(N)$.  Then
  $\ga_2\delta_2 \equiv \delta_2^\iota\ga_2' \pmod N$.  But by
  Shimura~\cite{shimura}, Lemma~3.29(1), we can then change $\ga_2$
  and $\ga_2'$ by elements of $\Ga_w(N)$ so that $\ga_2\delta_2 =
  \delta_2^\iota\ga_2'$, as desired.
\end{proof}

We can define a Petersson inner product on the space of weight
$(k_1,k_2)$ cusp forms just as in the one-variable case:
\[
\bracket{f,g} = \int_{\Gise(N)\backslash\uhh} f(z_i)\overline{g(z_i)}
y_1^{k_1-2}y_2^{k_2-2} \dx_1\dx_2\dy_1\dy_2
\]
(where $z_i = x_i+\sqrt{-1}y_i$); then just as in
Shimura~\cite{shimura}, Formula~(3.4.5), we see that the Hecke
operators $\Gise(N)(\delta_1,\delta_2)\Gise(N)$ and
$\Gise(N)(\delta_1^\iota,\delta_2^\iota)\Gise(N)$ are adjoint with
respect to that inner product.  Thus:

\begin{cor}
  \label{cor:hecke-alg-struct}
  The $\zz$-algebra generated by the Hecke operators is a commutative
  algebra; the Hecke operators are self-adjoint with respect to the
  Petersson inner product on $S_k(\Gise(N))$ and simultaneously
  diagonalizable.
\end{cor}

\begin{proof}
  The self-adjointness follows from Proposition~\ref{prop:self-adjoint}
  by the above discussion; the commutativity follows from
  Proposition~\ref{prop:self-adjoint} and Shimura~\cite{shimura},
  Proposition~3.8, and the simultaneous diagonalizability follows from 
  the self-adjointness.
\end{proof}

The effect of Hecke operators on Fourier expansions is given as
follows:

\begin{prop}
  \label{prop:hecke-expansion}
  Let $f$ be an element of $M_{(k_1,k_2)}(\Gise(N))$; if $a$ is an
  element of $\zmods{N}$, let $f|_{\smbigparens{\sigma_a,\smallmattwo
      1001}}$ have the Fourier expansion
  \[
  f|_\smbigparens{\sigma_a,\smallmattwo 1001}(z_1,z_2) =
  \sum_{m_1,m_2 \ge 0} c_{a,m_1,m_2}q_1^{m_1}q_2^{m_2}.
  \]
  If we set
  \[
  \Tni f(z_1,z_2) = \sum_{m_1,m_2 \ge 0} d_{m_1,m_2}q_1^{m_1}q_2^{m_2}
  \]
  then the $\dmi$'s are given by
  \[
  d_{m_1,m_2} = \sum_{\substack{a_1,a_2>0\\ a_i|(m_i,n_i)}}
  a_1^{k_1-1}a_2^{k_2-1} c_{(a_1/a_2),m_1n_1/a_1^2,m_2n_2/a_2^2}.
  \]
\end{prop}

\begin{proof}
  The proof is entirely parallel to the proof of the analogous fact in
  the one-variable case; c.f. Shimura~\cite{shimura}, (3.5.12).
\end{proof}

Note that the matrices $\smbigparens{\sigma_a,\smallmattwo 1001}$
don't normalize $\Gise(N)$.  This is why we have to introduce the
functions $f_a$ instead of simply diagonalizing $M_k(\Gise(N))$.

\begin{cor}
  \label{cor:eigenval-coeff}
  Let $f \in M_k(\Gise(N))$ be a simultaneous eigenform for all
  of the Hecke operators.  Then if $\lambda_{m_1,m_2}(f)$ is the
  eigenvalue for $T_{m_1,m_2}$, we have
  \[
  c_{m_1,m_2}(f) = \lambda_{m_1,m_2}(f)c_{1,1}(f).
  \]
  \qed
\end{cor}

Unfortunately, this Corollary isn't quite as useful as one might
hope, since the above coefficients are all zero by
Proposition~\ref{prop:zero-coeff} unless $\eps = -1$!  However, in
that situation, we do get the following result:

\begin{cor}
  \label{cor:hecke-eigenval}
  If $f$ and $g$ are elements of $S_k(\Gis{-1}(N))$ that are
  eigenfunctions for all $\Tni$'s with the same eigenvalues
  then, considered as elements of $\SkNb{-1}$, they differ by a
  multiplicative constant.
\end{cor}

\begin{proof}
  By Proposition~\ref{prop:zero-coeff} and
  Corollary~\ref{cor:eigenval-coeff}, if $c = c_{1,1}(f)/c_{1,1}(g)$
  then $c_{m_1,m_2}(f-cg)$ is zero unless $(m_i,N) > 1$.
\end{proof}

This can be restated as follows: let $\httkeb(N)$ be the $\cc$-algebra of
endomorphisms of $\SkeNb$ generated by the Hecke operators
$\Tni$ for $n_1 \equiv n_2 \pmod N$.  Then:

\begin{prop}
  \label{prop:cusp-free}
  The space $\SkNb{-1}$ is a free module of rank one over
  $\httkb{-1}(N)$.
\end{prop}

\begin{proof}
  By Corollary~\ref{cor:hecke-alg-struct}, we can find a basis for
  $\SkNb{-1}$ consisting of simultaneous eigenforms for all of the
  elements of $\httkb{-1}(N)$.  Furthermore, by
  Corollary~\ref{cor:hecke-eigenval}, no two of those eigenforms have
  the same eigenvalues.  This implies our Proposition.
\end{proof}

Similarly, we define $\httkes(N)$ to be the $\cc$-algebra of
endomorphisms of $\SkeN$ generated by the Hecke operators $\Tni$ for
$n_1 \equiv n_2 \pmod N$.  Proposition~\ref{prop:p-nobar} tells us
that the spaces $\Step$ and $\Stepb$ are equal; thus, the above
Proposition has the following Corollary:

\begin{cor}
  \label{cor:cusp-free}
  The space $\Stgisp{-1}$ is a free module of rank one over
  $\httsarg{(2,2)}{-1}(p)$.  \qed
\end{cor}

With a little bit more care, we can use the above techniques to prove
similar facts for $\eps = -k^2$ instead of just $\eps = -1$.  (This
isn't too surprising, since $\Xis{-1}(N)$ and $\Xis{-k^2}(N)$ are
isomorphic.)  They are in fact true for arbitrary $\eps$; the proof
demands different techniques, and will be given as
Theorem~\ref{thm:cusp-free}.  It does seem that $\Xis{-1}(N)$ is the
``dominant'' $\XiseN$; see Sections \ref{sec:cusp-relationships}
and~\ref{sec:hecke-kernel} for further discussion of this matter.

Finally, we let $\htteqs(N)$ denote the free polynomial algebra over
$\cc$ with variables $T_{n_1,n_2}$ for every pair $n_1$,$n_2$ of
positive integers that are relatively prime to $N$ and congruent mod
$N$.  This algebra acts on the spaces $\SkeN$ and $\SkeNb$ for all $k$
and $\eps$; its image in the endomorphism rings of those spaces gives
us the algebras $\httkes(N)$ and $\httkeb(N)$ that we defined above.

\section{Hecke Operators on $\XisN$}
\label{sec:hecke-two}

The Hecke operators $\Tni$ defined above have the following modular
interpretation: let $(E_1,E_2,\phi)$ be a point of $\XiseN$, and let
$\fdef{\pi_i}{E_i}{E_i'}$ be maps of elliptic curves of degree $n_i$,
where $(n_i,N) = 1$.  Then $\phi$ induces a map from $E_1'[N]$ to
$E_2'[N]$ which is an isomorphism of group schemes; $\Tni$ sends our
point to the sum of all points $(E_1', E_2', \phi)$ that arise in such
a fashion.  Why, then, do we impose the restriction that $n_1$ be
congruent to $n_2$ mod $N$?  The answer is that, if $\fdef{\pi}E{E'}$
is a map of degree $n$ (with $(n, N) = 1$) then $\pi$ doesn't preserve
the Weil pairing:
\begin{align*}
  (\pi x, \pi y) &= (x, \pi^\vee\pi y) \\
  &= (x, [n]y) \\
  &= (x,y)^n.
\end{align*}
So if $\phi$ raises the Weil pairing to the $\eps$'th power then, if
we push it forward via maps of order $n_i$ as above, the resulting map
raises the Weil pairing to the $\eps n_2/n_1$ power.  This
explains why we had to assume that $n_1 \equiv n_2 \pmod N$ for the Hecke
operators to act on the surfaces $\XiseN$.  However, we should have
Hecke operators $\Tni$ for arbitrary $n_i$ with $(n_i, N) = 1$ which
act on the surface $\XisN$.

The above considerations, when translated into matrices, lead us to
the following definition: for any $\eps$, $\eps'$ in $\zmods{N}$,
set
\[
\Disees(N) = \setof{\genpairmattwo}{
  \begin{array}{l}
    a_i,b_i,c_i,d_i \in \zz, \\
    a_id_i-b_ic_i > 0, \\
    (a_id_i-b_ic_i,N) = 1, \\
    \begin{array}{rcll}
      a_1 & \equiv & a_2 & \pmod N,\\
      b_1 & \equiv & \eps' b_2 & \pmod N, \\
      \eps c_1 & \equiv & c_2 & \pmod N, \\
      \eps d_1 & \equiv & \eps' d_2 & \pmod N
    \end{array}
  \end{array}
  }.
\]
It is obvious from the definitions that $\Diss{\eps,\eps} = \Dises$ and
one easily checks that
\[
\Disees\cdot\Diss{\eps',\eps''} \subset \Diss{\eps,\eps''}.
\]
These facts imply in particular that $\Disees$ is invariant under
multiplication by $\Gise(N)$ on the left and by $\Gis{\eps'}(N)$ on
the right; thus, $\Disees$ can be partitioned into Hecke operators
that send forms on $\Xise(N)$ to forms on $\Xis{\eps'}(N)$.  For any
$n_1$ and $n_2$ with $(n_i, N) = 1$ and with $\eps n_1 \equiv \eps'
n_2 \pmod N$, we define the Hecke operator $\Tni$ to be the sum of the
double cosets $\Gise(N)(\ga_1,\ga_2)\Gis{\eps'}(N)$ occurring in
$\Disees$ for which $\det(\ga_i) = n_i$.  This does depend on $\eps$,
but it has a set of left coset representatives that is independent of
$\eps$:

\begin{prop}
  \label{prop:hecke-2-coset}
  Let $n_1$ and $n_2$ be positive integers that are relatively prime
  to $N$, and let $\eps$ and $\eps'$ be elements of $\zmods{N}$ such
  that $\eps n_1 \equiv \eps' n_2 \pmod N$.  Then the set of elements
  of $\Disees(N)$ that have determinant $(n_1,n_2)$ has the following
  left coset decomposition:
  \[
  \coprod_{\substack{a_1,a_2>0\\a_id_i=n_i\\0\le b_i<d_i}}
  \Gise(N)\bigparens{\sigma_{a_1}\mattwo{a_1}{b_1N}{0}{d_1},
    \sigma_{a_2}\mattwo{a_2}{b_2N}{0}{d_2}}
  \]
  where, for $a \in \zmods{N}$, $\sigma_a$ is any matrix that is
  congruent to $\smallmattwo{a^{-1}}00a$ mod $N$.  Furthermore, the
  above left cosets are also disjoint as $\Ga(1)\times\Ga(1)$ cosets.
\end{prop}

\begin{proof}
  The proof is the same as the proof of
  Proposition~\ref{prop:hecke-coset}.
\end{proof}
Recall that we defined
\[
\SkisN = \prod_{\eps \in \zmods{N}} S_k(\Gise(N))
\]
and made a similar definition for $\SkisNb$.  Also, if $\fff$ is an
element of $\SkisN$, we write $\fff_\eps$ for its $\eps$'th component.
We then define Hecke operators $\Tni$ acting on the space $\SkisN$ by
setting $(\Tni \fff)_\eps = \Tni(\fff_{\eps n_2/n_1})$;
Proposition~\ref{prop:hecke-2-coset} shows that that action ``looks
the same'' for all $\eps$.  The following Proposition shows that the
action of these Hecke operators descends to the spaces $\SkeNb$, and
hence allows us to similarly define an action of them on the space
$\SkisNb$:

\begin{prop}
  If $f$ is a form in $S_k(\Gise(N))$ such that
  $c_{m_1,m_2}(f) = 0$ unless $(N, m_i) > 1$ then $\Tni f$ has
  the same property for all $n_i$ relatively prime to $N$.
\end{prop}

\begin{proof}
  The proof is the same as the proof of
  Proposition~\ref{prop:hecke-equiv}.
\end{proof}

The action on Fourier expansions is also as expected from
Proposition~\ref{prop:hecke-expansion}, with the same proof:

\begin{prop}
  \label{prop:hecke-2a-expansion}
  Let $f$ be an element of $M_{(k_1,k_2)}(\Gise(N))$; if $a$ is an
  element of $\zmods{N}$, let $f|_{\smbigparens{\sigma_a,\smallmattwo
      1001}}$ have the Fourier expansion
  \[
  f|_\smbigparens{\sigma_a,\smallmattwo 1001}(z_1,z_2) =
  \sum_{m_1,m_2 \ge 0} c_{a,m_1,m_2}q_1^{m_1}q_2^{m_2}.
  \]
  If we set
  \[
  \Tni f(z_1,z_2) = \sum_{m_1,m_2 \ge 0} d_{m_1,m_2}q_1^{m_1}q_2^{m_2}
  \]
  then the $\dmi$'s are given by
  \[
  d_{m_1,m_2} = \sum_{\substack{a_1,a_2>0\\a_i|(m_i,n_i)}}
  a_1^{k_1-1}a_2^{k_2-1} c_{(a_1/a_2),m_1n_1/a_1^2,m_2n_2/a_2^2}.
  \]
  \qed
\end{prop}

This Proposition (or Proposition~\ref{prop:hecke-2-coset}, which it is
a corollary of) allows us to translate theorems about forms on
$X_w(N)$ into theorems about forms on $\XisN$: if $f$ is a form on
some $\Xise(N)$ and we have a Hecke operator $\Tni$, we can consider
$f$ to be form on $X_w(N)\times X_w(N)$ and apply $T_{n_1} \times
T_{n_2}$ to it there.  This gives us a form on $X_w(N) \times X_w(N)$;
but by Proposition~\ref{prop:hecke-2-coset}, that has the same effect
as directly applying the $\Tni$ that we have defined above to $f$
considered as a form on $\Xise(N)$, so our resulting form, which is a
priori only a form on $X_w(N) \times X_w(N)$, is really a form on
$\Xis{\eps n_1/n_2}(N)$.  Thus, the fact that the Hecke operators
$T_n$ (with $(n,N) = 1$) on $X_w(N)$ commute implies that our Hecke
operators $\Tni$ commute.  Similarly, we can define a Petersson inner
product on $\SkisN$ by taking the orthogonal direct sum of the inner
products on the $\SkeN$'s; our Hecke operators are then normal with
respect to that inner product because the Hecke operators on $X_w(N)$
are.

It is frequently useful to encapsulate this relation between forms on
$\XisN$ and forms on $X_w(N)$ by defining a map
$\fdef{\Sigb}{\SkisNb}{\SkkGwNb}$ which sends $\fff \in \SkisNb$ to
$\sum_{\eps \in \zmods{N}} \fff_\eps$.  By $\SGwNb{k_i}$ we mean
$\SGwN{k_i}/V$ where $V$ is the space of forms $f \in \SGwN{k_i}$ such
that $c_m(f) = 0$ unless $(m,k_i) > 1$; it is a module over the Hecke
algebra generated by the operators $T_n$ with $(n,N) = 1$, and its
eigenspaces for that algebra are one-dimensional.  The following two
Propositions then sum up the discussion of the previous paragraph:

\begin{prop}
  \label{prop:xisn-xn}
  The map from $\SkisN$ to $\SkkGwN$ that sends a form $\fff$ to
  $\sum_{\eps\in \zmods{N}} \fff_\eps$ commutes with the action of
  Hecke operators.  It descends to an injection
%  $\fdef{\Sigb}{\SkisNb}{\SkkGwNb}$;  % formats badly
  $\Sigb$ from $\SkisNb$ to $\SkkGwNb$; if $\fff \in \SkisNb$ then
  \[
  \fff_\eps = \sum_{\substack{m_1,m_2 > 0\\\eps m_1+m_2\equiv 0
      \narpmod N\\(m_i,N)=1}} \cmi(\Sigb \fff)q_1^{m_1}q_2^{m_2}.
  \]
\end{prop}

\begin{proof}
  The only parts that remain to be proved are that $\Sigb$ is an
  injection and that $\fff_\eps$ can be recovered in the given manner.
  First, we note that, for all $m_1$, $m_2$ with $(m_i,N) = 1$,
  \[
  \cmi(\Sigb \fff) = \sum_{\eps \in \zmods{N}} \cmi(\fff_\eps).
  \]
  But Proposition~\ref{prop:zero-coeff} says that $\cmi(\fff_\eps) =
  0$ unless $\eps \equiv -m_2/m_1 \pmod N$; $\cmi(\Sigb \fff)$
  therefore equals $\cmi(\fff_{-m_2/m_1})$.  This together with
  Proposition~\ref{prop:zero-coeff} immediately implies our formula
  for $\fff_\eps$.  And if $\Sigb \fff = 0$ then this implies that,
  for all $\eps$ and for all $m_i$ such that $\eps \equiv -m_2/m_1
  \pmod N$, $\cmi(\fff_\eps)$ is zero.  But that implies that
  $\fff_\eps = 0$ by using Proposition~\ref{prop:zero-coeff} again.
\end{proof}

\begin{prop}
  \label{prop:hecke-alg-struct}
  The $\zz$-algebra generated by the Hecke operators $\Tni$ acting on
  $\SkisN$ is a commutative algebra; the Hecke operators are normal
  with respect to the Petersson inner product on $\SkisN$ and
  simultaneously diagonalizable.
\end{prop}

\begin{proof}
  This follows from the above reduction of these facts to facts about
  forms on $X_w(N)$ and from Shimura~\cite{shimura}, Theorem~3.41.
\end{proof}

Let $\fff$ be an element of $\SkisN$, and let $m_1$ and $m_2$ be
integers relatively prime to $N$.  We define $\cmi(f)$ to be equal to
$\cmi(\fff_{-m_2/m_1})$.  We also make the same definition for $\fff
\in \SkisNb$.  If we set $f = \sum_{\eps \in \zmods N} \fff_\eps$ then
$f$ is a form on $X_w(N) \times X_w(N)$, and $\cmi(\fff) = \cmi(f)$,
by Proposition~\ref{prop:zero-coeff}, as noted in the proof of
Proposition~\ref{prop:xisn-xn}.

\begin{prop}
  \label{prop:hecke-2-expansion}
  Let $\fff$ be an element of $\SkisN$; for $a \in \zmods{N}$, let
  $\fff_a$ be defined by
  \[
  (\fff_a)_\eps =
  \fff_{(a^{-2}\eps)}|_{\smbigparens{\sigma_a,\smallmattwo1001}}.
  \]
  Then for all $n_1$, $n_2$ with $(n_i,N) = 1$ and for all $m_1$,
  $m_2$ with $(m_i, N) = 1$, we have
  \[
  \cmi(\Tni\fff) = \sum_{\substack{a_1,a_2>0\\a_i|(m_i,n_i)}}
  a_1^{k_1-1}a_2^{k_2-1}
  c_{m_1n_1/a_1^2,m_2n_2/a_2^2}(\fff_{a_1/a_2}).
  \]
\end{prop}

\begin{proof}
  This is a corollary of Proposition~\ref{prop:hecke-2a-expansion}.
\end{proof}

We define $\htts(N)$ to be the free polynomial algebra over $\cc$ with
generators $\Tni$ for each pair $n_1$, $n_2$ of positive integers that
are relatively prime to $N$.  We define $\httkiss(N)$ to be its image
in the endomorphism ring of $\SkisN$; we define $\httkisb(N)$ to be
its image in the endomorphism ring of $\SkisNb$.

\begin{cor}
  \label{cor:big-eigenval-coeff}
  If $\fff \in \SkisN$ is a simultaneous eigenform for all Hecke
  operators $\Tni$ in $\httkiss(N)$ with eigenvalues $\lni(\fff)$
  then, for all $m_1$ and $m_2$ with $(m_i,N) = 1$, we have
  \[
  \cmi(\fff) = \lmi(\fff) c_{1,1}(\fff).
  \]
  \qed
\end{cor}

Thus, if $\fff$ is a non-zero element of $\SkisNb$ that is an
eigenform for all the $\Tni$'s then $c_{1,1}(\fff)$ is also non-zero;
we call such an $\fff$ a \emph{normalized eigenform} if $c_{1,1}(\fff)
= 1$.

\begin{cor}
  \label{cor:big-hecke-struct}
  The space $\SkisNb$ is a free module of rank one over $\httkisb(N)$.
\end{cor}

\begin{proof}
  By Proposition~\ref{prop:hecke-alg-struct}, we can find a basis for
  $\SkisNb$ consisting of simultaneous eigenforms for all elements of
  $\httkisb(N)$; the previous Corollary shows that the eigenspaces are
  one-dimensional, implying this Corollary.
\end{proof}

\begin{cor}
  \label{cor:big-hecke-struct-p}
  The space $\Sistp$ is a free module of rank one over
  $\httiss{(2,2)}(p)$.
\end{cor}

\begin{proof}
  This follows from Corollary~\ref{cor:big-hecke-struct} and
  Proposition~\ref{prop:p-nobar}.
\end{proof}

There is a special class of operators contained in our Hecke algebras
$\httkiss(N)$.  Given elements $\eps$ and $a$ of $\zmods N$, we have
\[
(1, \sigma_a)^{-1}\Gise(N)(1, \sigma_a) = \Gis{a^{-2}\eps}(N).
\]
The action of $(1, \sigma_a)$ therefore gives an isomorphism from
$\SkeN$ to $\SkN{a^{-2}\eps}$, denoted by $\diam{a}$; as with the
operators $\Tni$, $\diam{a}$ extends to the spaces $\SkisN$ and
$\SkisNb$ via the definition $(\diam{a}\fff)_\eps =
\diam{a}(\fff_{a^2\eps})$.  Furthermore, the action is the same if we
multiply $(1, \sigma_a)$ by $\pairsmallmattwo 1001a00a$; but if we
consider it as an operator on $X_w(N) \times X_w(N)$, as in the
discussion before Proposition~\ref{prop:xisn-xn}, then this, up to a
constant, is the product of the identity with the Hecke operator
$T(a,a)$.  By Shimura~\cite{shimura}, Theorem~3.24(4), $T(a,a)$ is in
the $\qq$-algebra generated by the $T(n)$'s, so $\diam{a}$ is in
$\httkiss(N)$.  Thus:

\begin{prop}
  \label{prop:diamond}
  For all $a \in \zmods N$, the operator $\diam{a}$ given by the
  action of $(1, \sigma_a)$ is an isomorphism from $\SkeN$ to
  $\SkN{a^{-2}\eps}$; furthermore, it is contained in $\httkiss(N)$.
  \qed
\end{prop}

\section{Relationships between the Spaces $\SkisNb$, $\SkeNb$, and
  $\SkNb{-1}$}
\label{sec:cusp-relationships}

When trying to prove that Hecke eigenspaces in $\SkeNb$ are
one-dimensional, we ran into problems because forms are ``missing''
Fourier coefficients: in particular, they don't have a $(1,1)$ Fourier
coefficient unless $\eps \equiv -1 \pmod N$, so we couldn't simply use
Corollary~\ref{cor:eigenval-coeff}.  However, the space $\SkisNb$
doesn't have that problem, and there is a natural projection map from
$\SkisNb$ to $\SkeNb$.  This gives us a replacement for the missing
Fourier coefficients; it also gives us a framework for seeing how the
spaces $\SkeNb$ differ (as $\htteqs(N)$-modules) as $\eps$ varies.

The key Lemma here is the following:

\begin{lemma}
  \label{lemma:big-to-small}
  The space $\SkeNb$ has a basis consisting of simultaneous
  $\httkeb(N)$-eigen\-forms $f$ that are of the form $\fff_\eps$ for
  simultaneous $\httkisb(N)$-eigenforms $\fff \in \SkisNb$.
\end{lemma}

\begin{proof}
  If $\fff \in \SkisNb$ is a $\httkisb(N)$-eigenform then it is
  certainly an eigenform for those Hecke operators $\Tni$ where $n_1
  \equiv n_2 \pmod N$; its $\eps$-component $\fff_\eps$ is therefore
  an eigenform for those operators as well.  The Lemma then follows
  from the fact that $\SkisNb$ has a basis of eigenforms, by
  Proposition~\ref{prop:hecke-alg-struct}.
\end{proof}

It is possible for two different $\httkisb(N)$-eigenforms in $\SkisNb$
to project to the same $\httkeb(N)$-eigenform in $\SkeNb$; we shall
discuss this in Theorem~\ref{thm:eigenform-lift}.  Also, some
eigenforms in $\SkisNb$ project to zero for some choices of $\eps$:
see the comments after the proof of the following Proposition and
Section~\ref{sec:hecke-kernel}.  We shall state a slightly stronger
version of this Lemma as Corollary~\ref{cor:big-to-small-uniquely}.

\begin{prop}
  \label{prop:eps-to-minus-one}
  If $f \in \SkeNb$ is a $\httkeb(N)$-eigenform then there is an
  $\httkb{-1}(N)$-eigenform $g \in \SkNb{-1}$ such that $\cmi(g) =
  \lmi(f)$ for all $m_1 \equiv m_2 \pmod N$.
\end{prop}

\begin{proof}
  By Lemma~\ref{lemma:big-to-small}, there is an eigenform $\fff \in
  \SkisNb$ such that $\lmi(\fff) = \lmi(f)$ for all $m_1 \equiv
  m_2 \pmod N$.  (We might a priori not be able to assume that
  $\fff_\eps = f$; however, $f$ is a linear combination of eigenforms
  projecting from $\SkisNb$, so those eigenforms must have the same
  eigenvalues as $f$.)  We can assume that $\fff$ is normalized.  We
  then set $g = \fff_{-1}$; it is a normalized eigenform contained in
  $\SkNb{-1}$, and $\lmi(g) = \lmi(\fff) = \lmi(f)$.  But
  Corollary~\ref{cor:eigenval-coeff} then tells us that $\cmi(g) =
  \lmi(f)$.
\end{proof}

Define $\KkeNbp$ to be the subspace of $\SkNb{-1}$ generated by
eigenforms whose eigenvalues are those of an eigenform in $\SkeNb$;
define $\KkeNb$ to be the subspace of $\SkNb{-1}$ generated by
eigenforms which do \emph{not} arise in such a fashion.  The Hecke
algebra $\httkeb(N)$ is isomorphic to the image of $\httkb{-1}(N)$ in
the endomorphism ring of $\KkeNbp$: both actions are diagonalizable,
so the rings are isomorphic iff the same eigenvalues occur, which is
the case by the definition of $\KkeNbp$ and by
Proposition~\ref{prop:eps-to-minus-one}.  In fact, the spaces
$\KkeNbp$ and $\SkeNb$ are isomorphic as $\htteqs(N)$-modules, because
the eigenspaces in $\SkeNb$ are one dimensional; we shall prove this
fact later as Theorem~\ref{thm:cusp-free}.  Thus, $\KkeNb$ measures
the difference between $\SkNb{-1}$ and $\SkeNb$; we shall study this
space in Section~\ref{sec:hecke-kernel}.

Since the proof of Proposition~\ref{prop:eps-to-minus-one} involved
lifting eigenforms in $\SkeNb$ to eigenforms in $\SkisNb$, we'd like
to see how ambiguous the choice of such a lifting is.  The following
Theorem answers that question:

\begin{theorem}
  \label{thm:eigenform-lift}
  Let $\fff$ be an eigenform in $\SkisNb$, and let $H \subset
  \zmods{N}$ be the set of $\eps$ such that $\fff_{-\eps} \ne 0$.
  Then:
  \begin{enumerate}
  \item $H$ is a subgroup of $\zmods{N}$.
  \item $H$ depends only on $\fff_{-1}$.
  \item Every element of $\zmods{N}/H$ has order one or two.
  \item If $\ggg$ is another eigenform in $\SkisNb$ then $\ggg_{-1} =
    \fff_{-1}$ if and only if there is a character $\chi$ on $H$ such
    that $\ggg_{-\eps} = \chi(\eps)\fff_{-\eps}$ for all $\eps \in H$.
  \end{enumerate}
\end{theorem}

First, we prove two Lemmas that we shall need during the proof of the
Theorem.

\begin{lemma}
  \label{lemma:coeffs-nd}
  Let $\fff$ be an eigenform in $\SkisNb$ and $\eps$ an element of
  $\zmods N$ such that $\fff_\eps \ne 0$.  For any positive integers
  $m_1$ and $m_2$ there exist positive integers $n_1$ and $n_2$ such
  that $\eps n_1 + n_2 \equiv 0$ (mod $N$), $(n_i,m_i) = 1$ for $i \in
  \{1,2\}$, and $\cni(\fff_\eps) \ne 0$.
\end{lemma}

\begin{proof}
  By Proposition~\ref{prop:xisn-xn}, $\Sigb \fff$ is an eigenform in
  $\SkkGwNb$.  Since the eigenspaces in $\SGwNb{k_i}$ are
  one-dimensional, there must exist $f_i \in \SGwNb{k_i}$ such
  that $\Sigb \fff = f_1 \otimes f_2$.

  For any $\eps' \in \zmods{N}$, set
  \[
  f_{i,\eps'} = \sum_{\substack{n > 0 \\n \equiv \eps' \narpmod N}}
  c_n(f_i)q^n.
  \]
  It is also an element of $\SGwNb{k_i}$.  (This follows easily from
  Shimura~\cite{shimura}, Proposition~3.64.)  Then
  \[
  \fff_\eps = \sum_{\eps' \in \zmods{N}} f_{1,\eps'}\otimes
  f_{2,-\eps\eps'},
  \]
  by Proposition~\ref{prop:xisn-xn}.

  Since $\fff_\eps \ne 0$, there exists $\eps' \in \zmods N$ such that 
  $f_{1,\eps'}$ and $f_{2,-\eps\eps'}$ are both nonzero.  By 
  Lang~\cite{lang-mod-forms}, Theorem~VIII.3.1, there exist
  $n_i$ such that $(n_i,Nm_i) = 1$ and that $c_{n_1}(f_{1,\eps'})$ and 
  $c_{n_2}(f_{2,-\eps\eps'})$ are both
  non-zero.  But Proposition~\ref{prop:xisn-xn} then implies that
  $\cni(\fff_\eps) \ne 0$, as desired.
\end{proof}

\begin{lemma}
  \label{lemma:forms-nonvanish}
  Let $\fff$ be an eigenform in $\SkisNb$ and $\eps$ an element of
  $\zmods N$ such that $\fff_\eps \ne 0$.  Then $\fff_{-\eps^j}$ is
  non-zero for all $j$.  In particular, $\fff_{-1/\eps}$ is non-zero.
\end{lemma}

\begin{proof}
  We can assume that $\fff$ is a normalized eigenform.  Since
  $\fff_{-\eps}$ is non-zero, there is some coefficient $\lambda =
  \cmi(\fff)$ that is non-zero, where $(m_i,N) = 1$ and $\eps m_1
  \equiv m_2 \pmod N$.  We therefore have $\Tmi(\fff) = \lambda \fff$, 
  by Corollary~\ref{cor:big-eigenval-coeff},
  so for all $\eps' \in \zmods N$,
  \begin{align*}
    \lambda \fff_{-\eps'} &= (\Tmi\fff)_{-\eps'} \\
    &= \Tmi(\fff_{-\eps' m_2/m_1}) \\
    &= \Tmi(\fff_{-\eps'\eps}).
  \end{align*}
  In particular, setting $\eps' = \eps^j$, we see that
  \[
  \lambda \fff_{-\eps^j} = \Tmi(\fff_{-\eps^{j+1}}),
  \]
  so if $\fff_{-\eps^j}$ is non-zero then, since $\lambda$ also is,
  $\fff_{-\eps^{j+1}}$ is as well, and we have our Lemma by induction.
\end{proof}

\begin{proof}[Proof of Theorem~\ref{thm:eigenform-lift}.]
  We can assume that $\fff$ is a normalized eigenform.  To show that
  $H$ is a subgroup, let $\eps_1$ and $\eps_2$ be elements of $H$.
  Thus, there exist $n_{1,i}$ and $n_{2,i}$ (for $i = 1,2$) such that
  $c_{n_{1,i},n_{2,i}}(\fff_{-\eps_i})$ is non-zero; by
  Lemma~\ref{lemma:coeffs-nd}, we can assume that $(n_{1,1},n_{1,2}) = 
  (n_{2,1},n_{2,2}) = 1$, and by Proposition~\ref{prop:zero-coeff},
  $\eps_in_{1,i} \equiv n_{2,i} \pmod N$.
  
  By Corollary~\ref{cor:big-eigenval-coeff},
  $c_{n_{1,i},n_{2,i}}(\fff) = \lambda_{n_{1,i},n_{2,i}}(\fff)$.  But
  \[
  \lambda_{n_{1,1}n_{1,2},n_{2,1}n_{2,2}}(\fff) =
  \lambda_{n_{1,1},n_{2,1}}(\fff)\lambda_{n_{1,2},n_{2,2}}(\fff),
  \]
  by our assumption that $(n_{i,1},n_{i,2}) = 1$, and is therefore
  non-zero, as is the corresponding Fourier coefficient of $\fff$.
  This is a Fourier coefficient of $\fff_{\eps}$ for
  \begin{align*}
    \eps &\equiv -(n_{2,1}n_{2,2}/n_{1,1}n_{1,2}) \\
    &\equiv -(n_{2,1}/n_{1,1})(n_{2,2}/n_{1,2}) \\
    &\equiv -\eps_1\eps_2.
  \end{align*}
  Thus, $\eps_1\eps_2 \in H$, so $H$ is a subgroup of $\zmods{N}$.
  
  To see that every element of $\zmods{N}/H$ has order one or two,
  pick $a \in \zmods{N}$ and let $\fff \in \SkisNb$ be an eigenform.
  Then $(\diam{a}\fff)_{-1} = \diam{a}(\fff_{-a^2})$.  Since
  $\diam{a}$ is an invertible operator contained in $\httkisb(N)$, by
  Proposition~\ref{prop:diamond}, the fact that $\fff_{-1} \ne 0$
  implies that $(\diam{a}\fff)_{-1} \ne 0$ as well, so so $\fff_{-a^2} 
  \ne 0$ and $a^2 \in H$.
  
  To show that $H$ depends only on $\fff_{-1}$, it's enough to prove
  the last part of the Theorem.  We shall prove that if $\ggg$ is an
  eigenform such that $\ggg_{-1} = \fff_{-1}$ then there is a
  character $\chi$ on $H$ such that $\ggg_{-\eps} =
  \chi(\eps)\fff_{-\eps}$; the converse (i.e. that $\ggg$'s
  constructed in that fashion are eigenforms) follows easily from the
  definitions.
  
  Thus, assume that we have normalized eigenforms $\fff$ and $\ggg$
  such that $\fff_{-1} = \ggg_{-1}$; let $\eps$ be an element of $H$,
  so $\fff_{-\eps} \ne 0$.  By Lemma~\ref{lemma:forms-nonvanish},
  $\fff_{-(1/\eps)}$ is also non-zero.  There then exist $m_1$ and
  $m_2$ relatively prime to $N$ such that $m_1 \equiv \eps m_2 \pmod
  N$ and $\cmi(\fff) \ne 0$.  Therefore, $\lmi(\fff)$ is also
  non-zero.  And
  \begin{align*}
    \lmi(\fff)\fff_{-\eps} &= (\Tmi\fff)_{-\eps} \\
    &= \Tmi(\fff_{-\eps m_2/m_1}) \\
    &= \Tmi(\fff_{-1}) \\
    &= \Tmi(\ggg_{-1}) \\
    &= \lmi(\ggg)\ggg_{-\eps}.
  \end{align*}
  Since $\lmi(\fff)$ and $\fff_{-\eps}$ are both non-zero, this implies 
  that $\lmi(\ggg)$ and $\ggg_{-\eps}$ are also both non-zero, and that 
  if we define $\chi(\eps) = \lmi(\fff)/\lmi(\ggg)$ (for any choice of 
  $m_i$ such that $m_1 \equiv \eps m_2 \pmod N$ and such that
  $\cmi(\fff_{-1/\eps}) \ne 0$) then $\ggg_{-\eps} =
  \chi(\eps)\fff_{-\eps}$, as desired.  We then only have to show that 
  $\chi$ is a character, not just a function; that follows by using
  the same arguments that we used to show that $H$ was a subgroup.
\end{proof}

We now have all the tools necessary to prove that the spaces $\SkeNb$
are free of rank one over $\httkeb(N)$ for all $\eps \in \zmods N$.

\begin{theorem}
  \label{thm:cusp-free}
  For all $\eps \in \zmods{N}$, all of the $\httkeb(N)$-eigenspaces in
  $\SkeNb$ are one-dimensional, and the space $\SkeNb$ is a free
  module of rank one over $\httkeb(N)$.
\end{theorem}

\begin{proof}
  Pick a $\httkeb(N)$-eigenspace in $\SkeNb$.  By
  Lemma~\ref{lemma:big-to-small}, it has a basis consisting of
  eigenforms of the form $\fff_\eps$ where $\fff$ is a normalized
  eigenform in $\SkisNb$.  Thus, we need to show that if $\fff$ and
  $\ggg$ are normalized eigenforms in $\SkisNb$ such that $\fff_\eps$
  and $\ggg_\eps$ are in the same eigenspace then $\fff_\eps$ and
  $\ggg_\eps$ are in fact constant multiples of each other.  However,
  $\lni(\fff_\eps) = \lni(\fff) = \cni(\fff)$, for all $n_1 \equiv n_2
  \pmod N$, so the fact that $\fff_\eps$ and $\ggg_\eps$ have the same
  eigenvalues simply means that $\fff_{-1}$ and $\ggg_{-1}$ are equal.
  Theorem~\ref{thm:eigenform-lift} then implies that $\fff_\eps$ and
  $\ggg_\eps$ are multiples of each other.  Thus, the eigenspaces are
  one-dimensional, and $\SkeNb$ is indeed a free $\httkeb(N)$-module
  of rank one.
\end{proof}

The basic idea behind the proof of Theorem~\ref{thm:cusp-free} is
that, if we have a form in $\SkeNb$, we can use
Lemma~\ref{lemma:big-to-small} to fill in the Fourier coefficients
that are forced to vanish by Proposition~\ref{prop:zero-coeff}.  Of
course, it's often easiest just to work with $\SkisNb$ and $\XisN$
directly.  As usual, we have the following Corollary:

\begin{cor}
  \label{cor:cusp-free-two}
  For all $\eps \in \zmods p$, the space $\Step$ is a free module of
  rank one over $\httsarg{(2,2)}{\eps}(p)$.
\end{cor}

\begin{proof}
  This follows from Theorem~\ref{thm:cusp-free} and
  Proposition~\ref{prop:p-nobar}.
\end{proof}

We also have the following slight strengthening of
Lemma~\ref{lemma:big-to-small}:

\begin{cor}
  \label{cor:big-to-small-uniquely}
  For every eigenform $f \in \SkeNb$ there exists an eigenform $\fff \in
  \SkisNb$ such that $\fff_\eps = f$.
\end{cor}

\begin{proof}
  By Lemma~\ref{lemma:big-to-small}, $\SkeNb$ has a basis consisting
  of such eigenforms.  Since the eigenspaces are one-dimensional,
  however, every eigenform must be a multiple of one of those basis
  elements.
\end{proof}

And, finally, we have the facts that $\KkeNbp$ and $\SkeNb$ are
isomorphic as $\htteqs(N)$-modules and a geometric consequence of that
fact:

\begin{cor}
  \label{cor:minus-one-split}
  For all $\eps \in \zmods N$, $\SkNb{-1}$ is isomorphic to $\KkeNb
  \oplus \SkeNb$ as a module over $\htteqs(N)$.
\end{cor}

\begin{proof}
  By definition, $\SkNb{-1} = \KkeNb \oplus \KkeNbp$.  But $\KkeNbp$
  is a $\htteqs(N)$-module that is a direct sum of one-dimensional
  spaces corresponding to the Hecke eigenvalues occurring in $\SkeNb$;
  the Corollary then follows from Theorem~\ref{thm:cusp-free}.
\end{proof}

\begin{cor}
  \label{cor:genus-maximized}
  If $N$ is a power of a prime then the geometric genus of (a
  desingularization of) $\Xise(N)$ is maximized when $\eps = -1$.
\end{cor}

\begin{proof}
  Corollary~\ref{cor:cusp-dim-coh} and
  Proposition~\ref{prop:bar-nobar} allow us to reduce this Corollary
  to showing that, for all $\eps$ and for all $M | N$, the dimension
  of $\Sbar_{(2,2)}(\Gis{-1}(M))$ is at least as large as the
  dimension of $\Sbar_{(2,2)}(\Gise(M))$.  This in turn follows
  directly from the above Corollary.
\end{proof}

This Corollary is in fact true for all $N \le 30$, as can be seen by
examining the tables in Kani and Schanz~\cite{kani-modular}.
Conjecture~\ref{conj:level-change} would imply this Corollary for all
natural numbers $N$, since in that case
Proposition~\ref{prop:bar-nobar} would be true for all $N$.

\section{The Hecke Kernel}
\label{sec:hecke-kernel}

In the previous Section, we saw that, for all $\eps \in \zmods N$, we
can write $\SkNb{-1}$ as $\KkeNb \oplus \SkeNb$.  Thus, the key to
understanding modular forms in all of the $\SkeNb$'s is to understand
the space $\SkNb{-1}$; once we have that, we then need to understand
its subspaces $\KkeNb$.  The goal of the present section is to study
those subspaces, which we call ``Hecke kernels''.  Note that
Corollary~\ref{cor:genus-maximized} gives us a geometric
interpretation of these spaces in some situations.

We first give the alternate following characterizations of forms in
$\KkeNb$:

\begin{prop}
  \label{prop:kernel-characterization}
  Let $f$ be an eigenform in $\SkNb{-1}$ and let $\eps$ be an element
  of $\zmods N$.  The following are equivalent:
  \begin{enumerate}
  \item $f$ is in $\KkeNb$.
  \item For any or all eigenforms $\fff \in \SkisNb$ such that
    $\fff_{-1} = f$, $\fff_\eps = 0$.
  \item For all $n_1$, $n_2$ such that $\eps n_1 + n_2 \equiv 0 \pmod
    N$, $\Tni f = 0$.
  \item For all $m_1$, $m_2$, $n_1$, and $n_2$ with $n_1m_1 \equiv
    n_2m_2 \pmod N$, $\eps n_1 + n_2 \equiv 0 \pmod N$, and $(n_i,m_i)
    = 1$ for $i \in \{1,2\}$, we have $c_{n_1m_1,n_2m_2}(f) = 0$.
  \end{enumerate}
\end{prop}

\begin{proof}
  We can assume $f$ is a normalized eigenform.  First we, show the
  equivalence between 1 and 2: let $\fff$ be an eigenform in $\SkisN$
  such that $\fff_{-1} = f$, which we can find by
  Corollary~\ref{cor:big-to-small-uniquely}.  By
  Theorem~\ref{thm:eigenform-lift}, $\fff_\eps$ only depends on the
  choice of $\fff$ up to a non-zero constant multiple.  If $\fff_\eps
  \ne 0$ then $\fff_\eps$ is an eigenform in $\SkeNb$ whose
  eigenvalues are the same as those of $f$, hence are the same as the
  Fourier coefficients of $f$, so $f$ isn't in $\KkeNb$.  Conversely,
  if $f$ isn't in $\KkeNb$ then there exists an eigenform $g \in
  \SkeNb$ whose eigenvalues are the Fourier coefficients of $f$.
  Corollary~\ref{cor:big-to-small-uniquely} allows us to pick an
  eigenform $\ggg \in \SkisNb$ such that $\ggg_\eps = g$; multiplying
  it (and $g$) by a constant factor, we can assume that $\ggg$ is a
  normalized eigenform.  Then $\ggg_\eps$ and $\ggg_{-1}$ have the
  same eigenvalues, so $\ggg_{-1}$ is a multiple of $f$, by our
  assumption on $g$; $\ggg$ therefore gives us an eigenform in
  $\SkisNb$ such that $\ggg_{-1} = f$ and $\ggg_\eps \ne 0$, as
  desired.  By Theorem~\ref{thm:eigenform-lift}, this is independent
  of the choice of $\ggg$, justifying our use of the phrase ``any or
  all''.
  
  Next we show that 2 and 3 are equivalent.  Thus, we are given
  normalized eigenforms $f \in \SkNb{-1}$ and $\fff \in \SkisNb$ such
  that $f = \fff_{-1}$ and we want to show that $\fff_\eps = 0$ iff,
  for all $n_1$ and $n_2$ such that $\eps n_1 + n_2 \equiv 0 \pmod N$,
  $\Tni f = 0$.  First assume that $\fff_\eps = 0$.  By
  Lemma~\ref{lemma:forms-nonvanish}, $\fff_{1/\eps} = 0$.  Then for
  all $n_i$ as above,
  \begin{align*}
    \Tni f &= \Tni(\fff_{-1}) \\
    &= (\Tni \fff)_{-n_1/n_2} \\
    &= (\Tni \fff)_{1/\eps} \\
    &= \lni(\fff)\fff_{1/\eps} \\
    &= 0.
  \end{align*}
  Conversely, if $\Tni f = 0$ for all $n_i$ with $\eps n_1 + n_2
  \equiv 0 \pmod N$ then the above series of equalities shows that
  $\lni(\fff)\fff_{1/\eps}$ is always zero, or equivalently (by
  Corollary~\ref{cor:big-eigenval-coeff}), $\cni(\fff)\fff_{1/\eps} =
  0$.  If $\fff_\eps \ne 0$ then there exist such $n_i$ such that
  $\cni(\fff) \ne 0$; thus, $\fff_{1/\eps} = 0$, so $\fff_\eps$ is
  zero after all, by Lemma~\ref{lemma:forms-nonvanish}.
  
  Next we show that 3 implies 4.  Assume that, for all $n_1$ and $n_2$
  with $\eps n_1 + n_2 \equiv 0$ (mod $N$), $\Tni f = 0$.  Then, for
  all $m_1$ and $m_2$ with $(m_i,n_i) = 1$, we have
  $T_{m_1n_1,m_2n_2}(f) = \Tmi(\Tni(f)) \\ = 0$, so in particular that
  is true for $m_i$ with $(m_i,n_i) = 1$ and with $m_1n_1 \equiv
  m_2n_2 \pmod N$.  But Corollary~\ref{cor:eigenval-coeff} then
  implies that $c_{m_1n_1,m_2n_2}(f) = 0$.
  
  Finally, we show that 4 implies 2, so let $f$ be a normalized
  eigenform such that all such coefficients $c_{m_1n_1,m_2n_2}(f)$ are
  zero, and let $\fff \in \SkisNb$ be a lift of $f$.  Assume that
  $\fff_\eps \ne 0$.  Thus, there exist $n_1$ and $n_2$ with
  $\cni(\fff) \ne 0$, or, equivalently, $\lni(\fff) \ne 0$.  Then for
  all $m_1$ and $m_2$ with $(m_i,n_i) = 1$ and with $m_1n_1 \equiv
  m_2n_2 \pmod N$, or equivalently $(1/\eps)m_1 + m_2 \equiv 0 \pmod
  N$,
  \begin{align*}
    0 &= \lambda_{m_1n_1,m_2n_2}(\fff) \\
    &= \lmi(\fff)\lni(\fff), \\
  \end{align*}
  so $\lmi(\fff) = 0$ for all $m_i$ with $(m_i,n_i) = 1$ and
  $(1/\eps)m_1 + m_2 \equiv 0 \pmod N$.  By
  Lemma~\ref{lemma:coeffs-nd}, $\fff_{1/\eps} = 0$; by
  Lemma~\ref{lemma:forms-nonvanish}, $\fff_\eps = 0$, a contradiction.
  Thus 4 implies 2.
\end{proof}

For an arbitrary form in $\KkeNb$, it is necessary for those
coefficients specified in part 4 of
Proposition~\ref{prop:kernel-characterization} to vanish.  The
following Proposition shows that even more coefficients of elements of
$\KkeNb$ vanish:

\begin{prop}
  \label{prop:kernel-diamond}
  For all $a$ and $\eps$ in $\zmods N$, the spaces $\KkeNb$ and
  $\KkNb{a^2\eps}$ are equal.
\end{prop}

\begin{proof}
  Let $f$ be an eigenform in $\KkeNb$; we want to show that $f$ is in
  $\KkNb{a^2\eps}$. Let $\fff$ be a lift of it to $\SkisNb$.  By
  Proposition~\ref{prop:kernel-characterization}, $\fff_\eps = 0$.
  Thus, $(\diam{a^{-1}}\fff)_{a^2\eps} = \diam{a^{-1}}(\fff_\eps)$ is
  also zero.  But by Proposition~\ref{prop:diamond}, $\diam{a^{-1}}$
  is in $\httkisb(N)$, so $(\diam{a^{-1}}\fff)$ is a multiple of
  $\fff$, which is non-zero since $\diam{a^{-1}}$ is invertible.  Thus,
  $\fff_{a^2\eps} = 0$, so $f$ is in $\KkNb{a^2\eps}$, by
  Proposition~\ref{prop:kernel-characterization}.
\end{proof}

Thus, if $\fff \in \SkisNb$ is a normalized eigenform such that
$\fff_{\eps}$ is zero for some $\eps$, or equivalently that
$\fff_{-1}$ is in $\KkeNb$, then $\fff_{a^2\eps}$ is also zero for all
$a \in \zmods N$.  So if we let $f = \Sigb \fff$ then lots of the
Fourier coefficients of $f$ are zero.  This leads one to suspect that
$f$ might be related to forms with complex multiplication, where we
define an eigenform $g$ on $X_w(N)$ to have \emph{complex
  multiplication} if there exists a non-trivial character $\phi$ such
that $\phi(p)\lambda_p(g) = \lambda_p(g)$ (or, equivalently,
$\lambda_p(g) = 0$ unless $\phi(p) = 1$) for all primes $p$ in a set
of density one, where $\lambda_p(g)$ is the $T_p$-eigenvalue for $g$.
(This is as in Ribet~\cite{ribet-nebentypus}, \S3, except that we
don't require $g$ to be a newform.)  We also say that $g$ is a
\emph{CM-form}.  It is indeed the case that such forms are linked to
elements of the Hecke kernel:

\begin{theorem}
  \label{thm:kernel-cm}
  An eigenform $f$ is in $\KeNb{(k_1,k_2)}$ if and only if there exist
  eigenforms $f_i \in S_{k_i}(\Ga_w(N))$ such that, for all $n_1 \equiv
  n_2 \pmod N$ with $(n_i,N) = 1$,
  \[
  \cni(f) = c_{n_1}(f_1)c_{n_2}(f_2)
  \]
  and such that the $f_i$ have complex multiplication by some
  character $\phi$ such that $\phi(-\eps) = -1$.  Furthermore,
  $\KeNb{(k_1,k_2)}$ is spanned by such forms.
\end{theorem}

\begin{proof}
  Let $k = (k_1,k_2)$, and let $f \in \SkeNb$ be an eigenform.  Pick
  an eigenform $\fff \in \SkisNb$ such that $\fff_{-1} = f$ and let
  $H$ be the subgroup of $\eps' \in \zmods N$ such that $\fff_{-\eps'}
  \ne 0$, as in Theorem~\ref{thm:eigenform-lift}.  By
  Proposition~\ref{prop:xisn-xn}, $\Sigb \fff$ is an eigenform in
  $\SkkGwNb$; but eigenspaces in that latter space are
  one-dimensional, so $\Sigb \fff = f_1 \otimes f_2$, where $f_i \in
  \SGwNb{k_1}$ is an eigenform.  We wish to relate $f$'s being an
  element of $\KkeNb$, i.e. having $\fff_{\eps} = 0$, to the $f_i$'s
  being CM-forms.
  
  For all $m_1$ and $m_2$ with $(m_i,N) = 1$, $\cmi(\fff) =
  c_{m_1}(f_1)c_{m_2}(f_2)$.  If $\eps' \nin H$, i.e. $\fff_{-\eps'} =
  0$, then, for all $m_i$ such that $\eps' m_1 \equiv m_2 \pmod N$,
  $\cmi(\fff) = 0$, so $c_{m_1}(f_1) = 0$ or $c_{m_2}(f_2) = 0$.
  Since the $f_i$ are eigenforms, their first Fourier coefficients are
  non-zero; thus, setting $m_2 = 1$, $c_{m_1}(f_1) = 0$ for $m_1
  \equiv 1/\eps' \pmod N$ where $\eps' \nin H$.  Since $H$ is a
  subgroup, this means that $c_{m_1}(f_1) = 0$ for $m_1 \nin H$
  (identifying $m_1$ with its projection to an element of $\zmods N$).
  Similarly, $c_{m_2}(f_2) = 0$ for $m_2 \nin H$.
  
  First, assume that $f \in \KkeNb$, i.e. that $\fff_\eps = 0$, or
  that $-\eps \nin H$.  Pick a non-trivial character $\phi$ of $\zmods
  N$ that is trivial on $H$ and such that $\phi(\eps) \ne -1$.  The
  previous paragraph shows that $f_1$ and $f_2$ both have complex
  multiplication by $\phi$.  By part 3 of
  Theorem~\ref{thm:eigenform-lift}, $\phi$ has order two; thus,
  $\phi(-\eps) = -1$, as desired.
  
  Conversely, assume that there exists a character $\phi$ such that
  the forms $f_i$ have complex multiplication by $\phi$ and such that
  $\phi(-\eps) = -1$.  Pick $m_1$ and $m_2$ such that $\eps m_1 + m_2
  \equiv 0$ (mod $N$).  Then $-\eps \equiv m_2/m_1 \pmod N$; since
  $\phi(-\eps) = -1$, either $\phi(m_1)$ or $\phi(m_2)$ is not equal
  to one.  Thus, either $c_{m_1}(f_1)$ or $c_{m_2}(f_2)$ is zero, so
  $\cmi(\fff) = 0$.  This is true for all such $m_i$, so $\fff_\eps =
  0$, i.e. $f \in \KkeNb$.

  Finally, the fact that $\KkeNb$ is spanned by such forms follows
  from the fact that it has a basis of eigenforms, which is obvious
  from the definition of $\KkeNb$.
\end{proof}

For $p$ prime we define $\Kis(p)$ to be the subspace
$\overline{K}_{(2,2),\eps}(p)$ of $\Step$ for any $\eps \in \zmods p$
such that $-\eps$ is non-square, where we identify $\Step$ with
$\Stepb$ by Proposition~\ref{prop:p-nobar}.  (For this to make sense,
we should assume that $p \ne 2$; since $\Stgis{\eps}{2}$ is zero for
all $\eps$, this isn't very important.)  This is independent of the
choice of $\eps$ by Proposition~\ref{prop:kernel-diamond}; its
dimension is the difference between the geometric genera of
$\Xis{-1}(p)$ and $\Xise(p)$, by Corollary~\ref{cor:genus-maximized}.
We shall give an explicit basis for this space in Sections
\ref{sec:prime} and~\ref{sec:examples}.

\section{The Adelic Point of View}
\label{sec:adelic}

As we have seen in Section~\ref{sec:hecke-two}, to get a satisfactory
theory of Hecke operators, we had to consider the surface $\XisN$, not
just the surfaces $\XiseN$.  To explain this, it helps to look at
$\XisN$ from the adelic point of view.  Thus, we review some of
definitions from that theory and explain their relevance to our
context.  For references, see Diamond and Im~\cite{diamond-im},
Section~11.

Let $\ado$ denote the finite adeles, i.e. the restricted direct
product of the fields $\qq_p$ with respect to the rings $\zz_p$.  Let
$U$ be an open compact subgroup of $\GLtA$.  We define the curve $Y_U$
to be $\GLtpQ \backslash (\uhG) / U$.  Here, $\GLtpQ$ is the set of
matrices in $\GL_2(\qq)$ with positive determinant, acting on $\uh$
via fractional linear translations and on $\GLtA$ via the injection
$\qq \into \ad$; $U$ acts trivially on $\uh$ and acts on $\GLtA$ via
multiplication on the right.  This defines $Y_U$ as a non-compact
curve over the complex numbers; it has a canonical compactification
$X_U$ given by adding a finite number of cusps.  The curves $X_U$ and
$Y_U$ in fact have canonical models over $\qq$ which are irreducible;
over $\cc$, however, the number of their components is given by the
index of $\det U$ in $\zzhat^\times$.  If $U$ and $U'$ are open
compact subgroups of $\GLtA$ and if $g$ is an element of $\GLtA$ such
that $g^{-1}Ug \subset U'$ then multiplication by $g$ on the right
gives a map $\fdef{g^*}{X_U}{X_U'}$; it descends to the models over
$\qq$.

We define a \emph{cusp form of weight $k$ on $X_U$} to be a function
$\fdef{\fff}{\uhG}{\cc}$ such that
\begin{enumerate}
\item $\fff(z,g)$ is a holomorphic function in $z$ for fixed $g$.
\item $\fff(\gamma z,\gamma g) = j(\gamma,z)^k\fff(z,g)$ for all
  $\gamma\in\GLtpQ$.
\item $\fff(z,gu) = \fff(z,g)$ for all $u \in U$.
\item $\fff(z,g)$, considered as a function in $z$, vanishes at
  infinity for all $g$.
\end{enumerate}
We denote by $S_k(U)$ the space of all such forms.  If $g^{-1}Ug
\subset U'$ then we get a map $\fdef{g_*}{S_k(U')}{S_k(U)}$ by
defining $(g_*\fff)(z,h)$ to be $\fff(z,hg)$.

Each $U$-double coset in $\GLtA$ gives a Hecke operator, which acts on
$S_k(U)$.  If $U = \GL_2(\zz_p) \times U^p$ then the Hecke operator
$T_p$ is generated by the elements of $\Mt(\zz_p)$ whose determinant
is in $p\zz_p^\times$; defining the Hecke operator $S_p$ to be the
double coset generated by $\smallmattwo p00p$ in the $\GLtQp$
component, the ring of Hecke operators consisting of those double
cosets generated by elements in $\GL_2(\qq_p)$ is generated by $T_p$
and $S_p^{\pm 1}$.

If we define $S_k(\cc)$ to be the direct limit of the $S_k(U)$'s as
$U$ gets arbitrarily small then the above maps $g_*$ make this into an
admissible representation of $\GLtA$; the original spaces $S_k(U)$ can
be recovered from that representation by taking its $U$-invariants.
The main fact that we need is the following adelic analogue of parts
of Atkin-Lehner theory:

\begin{theorem}[Strong Multiplicity One]
  If $\pi$ and $\pi'$ are two irreducible constituents of $S_k(\cc)$
  such that $\pi_p$ and $\pi'_p$ are isomorphic for almost all $p$
  then $\pi$ and $\pi'$ are equal.  (Not just isomorphic.)
  Furthermore, if $\fff$ and $\fff'$ are elements of $\pi$ and $\pi'$
  then this is the case iff $\fff$ and $\fff'$ have the same
  eigenvalues for almost all $T_p$ and $S_p$; in this case, they have
  the same eigenvalues for all $p$ such that $\fff \in S_k(U)$ for
  some $U$ of the form $\GL_2(\zz_p) \times U^p$.  \qed
\end{theorem}

The subgroups that we shall be concerned with are
\[
U_w(N) = \setof{g \in \GLtZh}{g \equiv \mattwo *001 \pmod N}
\]
and
\[
U(N) = \setof{g \in \GLtZh}{g \equiv \mattwo 1001 \pmod N}.
\]
These define the modular curves $X_w(N)$ and $X(N)$, respectively.
The modular interpretation of $X(N)$ is given as follows: for each
$\eps \in \zmods N$, choose a matrix $g_\eps \in \GLtZh$ congruent to
$\smallmattwo {\eps^{-1}}001$ mod $N$.  The strong approximation
theorem for $\GL_2$ implies that every point in $Y(N)$ has a
representative of the form $(z, g_\eps)$ for some unique choice of
$\eps$; we let this point correspond to the elliptic curve
$\cc/\bracket{z,1}$ together with the basis for its $N$-torsion given
by $(\eps z/N, 1/N)$.  We then have an action of $\GLtN$ on $X(N)$
that sends a matrix $\overline{g} \in \GLtN$ to the map
$\fdef{(g^{-1})^*}{X(N)}{X(N)}$, where $g$ is any lifting of
$\overline g$ to $\GLtZh$; it has the modular interpretation of
preserving the elliptic curve and having $\overline g$ act on the
basis for its $N$-torsion on the left.

Note that, in contrast, the action of $\SLtN$ on $X_w(N)$ can't easily
be defined adelically; this is one reason why one can't define such an
action over $\qq$, and thus why we find it convenient to use the
curves $X(N)$ rather than $X_w(N)$ at times.  However, with a bit of
care it is possible to use the action of $\GLtN$ on $X(N)$ to extract
information about the action of $\SLtN$ on $X_w(N)$; we shall do this
in Section~\ref{sec:prime}.

Now we turn to the surfaces $\XisN$.  Definitions similar to the above
go through, replacing $\uhG$ by $\uhhGG$ and putting in two copies of
everything else.  We then recover our surfaces $\XisN$ and spaces
$\SkisN$ of cusp forms by using the following subgroup:
\[
\UisN = \setof{(g_1,g_2) \in \GLtZZh}{g_1\equiv g_2 \pmod N}.
\]
The above definitions of Hecke operators pass over immediately to our
situation; in particular, it is easy to check that $T_{p_1,p_2}$ is
$T_{p_1} \times T_{p_2}$ (for $(p,N) = 1$) and $\diam {p}$ is $1
\times S_p$ (again for $(p,N) = 1$; note that $S_p \times 1$ is
$\diam{p^{-1}}$).  Using these definitions, we also easily see that
that, as claimed,
\[
\XisN = \GLtN\backslash(X(N)\times X(N)),
\]
where $\GLtN$ acts diagonally with the action given above.

In contrast with this situation, there does \emph{not} exist a
subgroup $U_{\simeq,\eps}(N)$ that would allow us to define $\XiseN$
in the same way; this explains why we couldn't naturally define a
Hecke operator $\Tni$ acting on $\XiseN$ unless $n_1 \equiv n_2 \pmod
N$.  Of course, it isn't hard to see which points on $\XisN$ are on
$\XiseN$ for some $\eps$: they are the points that have a
representative of the form $(z_1,z_2,g_1,g_2)$ with $g_i \in \GLtZh$
and with $\det g_1 \equiv \eps \det g_2 \pmod N$.  And if we are given
$\fff \in S_k(\UisN) = \SkisN$, we can recover $\fff_\eps$ from it by
letting
\[
\fff_\eps(z_1,z_2) = \fff(z_1,z_2,1,g_\eps).
\]

\section{The Case of Prime Level}
\label{sec:prime}

In this Section, we discuss facts that are special to the case of
weight $(2,2)$ forms on prime level.  The main fact here is that we
can ignore Fourier coefficients that are multiples of $p$, as stated
in Proposition~\ref{prop:p-nobar}; this in turn implies that certain
spaces of cusp forms are free of rank one over their Hecke algebras,
as stated in Corollaries \ref{cor:big-hecke-struct-p}
and~\ref{cor:cusp-free-two}.  In the rest of this Section, we shall
present some general calculations that lead us towards methods for
calculating the spaces $\Step$; in the next Section, we shall give
some explicit constructions of forms.

Since
\[
\Step = (\StGwp\otimes\StGwp)^\SLtp,
\]
to understand $\Step$ we should understand the representation theory
of $\SLtp$ on $\StGwp$.  Since $\smallmattwo {-1}00{-1}$ acts
trivially on $\StGwp$, we can look at the representation theory of
$\PSLtp$ instead.  We shall start by considering arbitrary weights and
levels, and adding the assumptions of weight 2 and level $p$ as it
becomes convenient.

The basic fact about representations of groups on spaces of cusp forms
is the Strong Multiplicity One Theorem.  This tells us how to pick out
the irreducible representations of $\GLtA$ that are contained in
$S_k(\cc)$: they are just the Hecke eigenspaces.  Taking
$\GLtN$-invariants, this breaks up $\SkUN$ into smaller
subrepresentations of $\GLtN$.  (Of course, these smaller
subrepresentations may not be irreducible as representations of
$\GLtN$.)  To apply this, we need to relate $\SkUN$ and its
eigenspaces to spaces that we understand better.

First we recall that ${\mattwo N001}^{-1}\Ga_w(N){\mattwo N001}
\subset \Ga_1(N^2)$.  This allows us to pass from forms on $X_w(N)$ to
forms on $X_1(N^2)$: the image of $\SkGwN$ is the direct sum of the
spaces $S_k(\Ga_0(N^2),\chi)$ where $\chi$ is a character on $\zmods
N$.  A form $f = \sum c_m q^m$, where $q = \expn{z/N}$, gets sent to a
form with the same Fourier expansion except that $q$ is now equal to
$\expn z$.  Furthermore, if $\psi$ is a character on $\zmods N$ then
the form $f_\psi$, which is defined to have Fourier expansion $\sum
c_m\psi(m)q^m$, is still a form in $\SkGwN$, by
Shimura~\cite{shimura}, Proposition~3.64.

We now turn to producing forms contained in $\SkUN$.  A form $\fff \in
\SkUN$ is a function from $\uhG$ to $\cc$ with those properties listed
in Section~\ref{sec:adelic}; it then follows easily that if, for $\eps
\in \zmods N$, we define $\fff_\eps$ by setting $\fff_\eps(z) =
\fff(z,g_\eps)$ (where $g_\eps$ is a matrix in $\GLtZh$ that is
congruent to $\smallmattwo {\eps^{-1}}001$ mod $N$) then each of the
$\fff_\eps$'s is a form in $\SkGwN$.  By the Strong Approximation
Theorem, a choice of such $\fff_\eps$'s determines $\fff$ uniquely.
Thus, we can think of forms on $\SkUN$ as $\phi(N)$-tuples of forms on
$\SkGwN$.

This allows us to determine the Hecke eigenspaces in $\SkUN$.  The
dimension of $\SkUN$ is $\phi(N)$ times the dimension of $\SkGwN$, so
the hope is that each eigenform on $\SkGwN$ will somehow give us
$\phi(N)$ different eigenforms on $\SkUN$.  This is indeed what
happens, as we shall see in Proposition~\ref{prop:skun-eigenforms}:

\begin{lemma}
  \label{lemma:un-gwn}
  Let $\fff$ be an element of $\SkUN$ and let $q$ be a prime not
  dividing $N$.  Then, for all $\eps \in \zmods N$, $(T_q \fff)_\eps =
  T_q(\fff_{\eps q})$ and $(S_q \fff)_\eps = S_q(\fff_{\eps q^2})$.
\end{lemma}

\begin{proof}
  This follows from tracing through the definitions; alternately one
  can use the modular interpretation of points on $X(N)$ and Hecke
  operators together with the fact that if $\fdef{\pi}{E}{E'}$ is an
  isogeny of degree $N$ then $(\pi x, \pi y)_{E'} = (x,y)_E^n$, where
  $(,)_E$ denotes the Weil pairing.
\end{proof}

\begin{cor}
  \label{cor:un-gwn}
  Let $g \in \SkGwN$ be an eigenform, with eigenvalues
  $\{a_q,\chi(q)\}$ (for $T_q$ and $S_q$ respectively, as $q$ varies
  over primes not dividing $N$).  Let $\psi$ be a character of $\zmods
  N$.  Then the form $\fff(g,\psi) \in \SkUN$ defined by
  $\fff(g,\psi)_\eps = \psi(\eps)g$ is an eigenform with eigenvalues
  $\{\psi(q)a_q, \psi^2(q)\chi(q)\}$.
\end{cor}

\begin{proof}
  Write $\fff$ for $\fff(g,\psi)$.  By the Lemma,
  \begin{align*}
    (T_q \fff)_\eps &= T_q(\fff_{\eps q}) \\
    &= T_q(\psi(\eps q)g) \\
    &= \psi(q)\psi(\eps)a_q g \\
    &= \psi(q)a_q \fff_\eps.    
  \end{align*}
  The calculation for $S_q$ proceeds in exactly the same manner.
\end{proof}

This allows us to produce a basis of eigenforms for $\SkUN$ in terms
of a basis of eigenforms for $\SkGwN$:

\begin{prop}
  \label{prop:skun-eigenforms}
  Let $\{g_j\}$ be a basis of eigenforms for $\SkGwN$.  Then the set
  of forms $\{\fff(g_j,\psi)\}$, as $g_j$ varies over elements of the
  basis and $\psi$ varies over characters of $\zmods N$, give a basis
  of eigenforms for $\SkUN$.  Every set $\{a_q,\chi(q)\}$ of
  eigenvalues for $T_q$ and $S_q$ (as $q$ runs over primes not
  dividing $N$) that occurs in $\SkUN$ occurs in $\SkGwN$.  A basis
  for the set of eigenforms in $\SkUN$ with eigenvalues $\{a_q,
  \chi(q)\}$ is given by taking the forms $\fff(g,\psi)$ where $\psi$
  varies over the characters of $\zmods N$ and where, once $\psi$ is
  fixed, $g$ varies over a basis for those eigenforms in $\SkGwN$
  which have eigenvalues $\{a_q\psi^{-1}(q), \chi(q)\psi^{-2}(q)\}$.
\end{prop}

\begin{proof}
  Assume that we have an expression of linear dependence involving the
  forms $\fff(g_j,\psi)$.  Looking at the first coordinate, the fact
  that the forms $\{g_j\}$ form a basis for $\SkGwN$ implies that we
  can assume that our relation involves only forms $\fff(g,\psi)$ for
  some fixed form $g$.  But those forms are linearly independent since
  characters are linearly independent.  This gives us $\phi(N)\cdot
  \dim \SkGwN$ forms; but that's the dimension of $\SkUN$, so those
  forms give a basis for $\SkUN$ that consists of eigenforms.

  Every set of eigenvalues on $\SkUN$ is therefore of the form
  $\{\psi(q)a_q, \psi^2(q)\chi(q)\}$, where $\{ a_q,\chi(q)\}$ is the
  set of eigenvalues of a form $g \in \SkGwN$, by
  Corollary~\ref{cor:un-gwn}.  But those are the eigenvalues of
  $g_\psi$, which is also an eigenform in $\SkGwN$.  The last 
  statement of the Proposition follows in a similarly direct manner
  from the first paragraph of the proof and
  Corollary~\ref{cor:un-gwn}.
\end{proof}

To restate the last sentence of the above Proposition: assume that $g
\in \SkGwN$ is a newform with eigenvalues $\{a_p,\chi(p)\}$.  A basis
for the eigenforms in $\SkUN$ with those eigenvalues is given by the
forms $\fff(g_{\psi^{-1}},\psi)$ together with the forms
$\fff(h,\psi)$ where $h$ runs over oldforms with the same eigenvalues
as $g_{\psi^{-1}}$.
  
Let us now fix $k=2$ and $N=p$ prime.  We may assume that $p > 5$,
since $\StGwp$ is zero otherwise.  Pick a set $A = \{a_q,\chi(q)\}$ of
eigenvalues.  Let $g \in \StGwp$ be a newform with those eigenvalues;
we wish to calculate the dimension of the space $S_A$ of forms in
$\StUp$ with eigenvalues $A$.  For each character $\psi$, we can
produce an element of $S_A$ all of whose components are multiples of
$g_{\psi^{-1}}$; this gives us $(p-1)$ forms.  Furthermore, when
$g_{\psi^{-1}}$ is an oldform, we can produce extra forms.  Since
$S_2(\Ga(1))$ is zero, we can produce at most one extra form for each
$\psi$ this way: this happens when the eigenvalues
$\{a_q\psi^{-1}(q),\chi(q)\psi^{-2}(q)\}$ occur in $S_2(\Ga_1(p))$.

For how many $\psi$ does an extra form arise in this way?  By the
Strong Multiplicity one theorem, studying $S_A$ reduces to the study
of irreducible representations of $\GLtA$ and their $U(p)$-invariants.
Factoring those representations, we have to study irreducible
representations of $\GL_2(\qq_q)$ and their $U(p)_q$-invariants.  If
$q \ne p$ then $U(p)_q = \GL_2(\zz_q)$; since the space of
$\GL_2(\zz_q)$ invariants of an irreducible representation of
$\GL_2(\qq_q)$ is either zero- or one-dimensional, we can therefore
concentrate on the irreducible representations of $\GL_2(\qq_p)$, and
in particular calculating the dimension of their $U(p)_p$-invariants,
where
\[
U(p)_p = \setof{g \in \GL_2(\zz_p)}{g \equiv \mattwo 1001 \pmod p}.
\]
Irreducible representations of $\GL_2(\qq_q)$ can be classified as
\emph{principal series}, \emph{special}, or \emph{supercuspidal}.  If
the space of $U(p)_p$-invariants is nonzero then it is $(p+1)$-, $p$-,
or $(p-1)$-dimensional, depending on which classification it falls
into; thus, we have two, one, or no extra dimensions of oldforms
arising in the principal series, special, and supercuspidal cases,
respectively.

Let us now turn towards the space $\StGwp$.  The group $\PSLtp$ acts
on this space; we wish to determine its irreducible representations.
Since this action is not given adelically, we can't just apply the
theory of irreducible $\GLtA$-representations and the Strong
Multiplicity One Theorem to get the answer.  However, we can use the
adelic action to get information about this representation as follows:
let $g$ be an element of $\StGwp$ and let $\fff$ be an element of
$\StUp$ such that $\fff_1 = g$.  Let $\overline{\ga}$ be an element of
$\PSLtp$ and let $\ga$ be an element of $\GLtZp$ projecting to it.
Then $\overline{\ga}$ sends $g$ to $(\ga^{-1}_* g)_1$, as can be seen
by tracing through the definitions.  In particular, we get
representations of $\PSLtp$ on $\StGwp$ by projecting the
representations given in the previous paragraphs down to their first
coordinate.

The map from $\StUp$ to $\StGwp$ sending $\fff$ to $\fff_1$ is
injective unless there is a $\psi$ such that $g = g_{\psi}$, by
Proposition~\ref{prop:skun-eigenforms}, i.e. unless $g$ is a CM-form,
in which case all of the forms in the representation are CM-forms, and
the dimension of the representation in $\StGwp$ is half of the
dimension of the representation in $\StUp$.  Thus, we have decomposed
$\StGwp$ as a direct sum of representations that are either of
dimension $p-1$, $p$, $p+1$, $(p-1)/2$, or $(p+1)/2$.

These representations may not be irreducible, however.  Most of the
time, they do turn out to be irreducible; we can see this by looking
at the character table of $\PSLtp$.  The dimensions of the irreducible
representations of $\PSLtp$ are 1, $p-1$, $p$, $p+1$, and either
$(p-1)/2$ (if $p \equiv 3 \pmod 4$) or $(p+1)/2$ (if $p \equiv 1 \pmod
4$).  Furthermore, the only one-dimensional representation of $\PSLtp$
is the trivial one, which doesn't occur in $\StGwp$ (since that would
be equivalent to having a form that is invariant under $\PSLtp$, i.e.
a form in $S_2(\Ga(1))$).  There are no $2$-dimensional
representations, either, so by comparing dimensions, we see that the
representations that we have constructed above are either trivial or
the direct sum of two representations of dimension $(p-1)/2$ or
$(p+1)/2$.

We wish to see how $\dim \Step$ varies as a function of $\epsilon$.
Write $\chi_w(p)$ for the character of $\StGwp$, considered as a
representation of $\PSLtp$.  Then
\begin{align*}
  \dim \Step &= \dim (\StGwp \otimes \StGwp\circ\theta_\epsilon)^\PSLtp \\
  &= \bracket{\chi_w(p) \otimes \chi_w(p) \circ \theta_\epsilon, 1} \\
  &= \bracket{\chi_w(p), \overline{\chi_w(p)} \circ \theta_\epsilon}.
\end{align*}
Assume that $\StGwp$ has $\bigoplus_i R_i^{\oplus n_i}$ as its
decomposition into a sum of irreducible representations.  Then, by the 
above,
\[
\dim \Step = \sum_{\substack{i,j\\ R_i \simeq \overline{R_j}\circ
  \theta_\epsilon}} n_in_j.
\]

Now assume that $p \equiv 1$ (mod $4$).  Examining the character table
of $\PSLtp$, we see that $R_i \simeq \overline{R_i}$ for all $R_i$ and
that $R_i \simeq R_i \circ \theta_\epsilon$ for all $\epsilon$ unless
$R_i \simeq W'$ or $W''$, where $W'$ and $W''$ are the irreducible
representations of dimension $(p+1)/2$.  In this latter case,
composing with $\theta_\epsilon$ switches $W'$ and $W''$ if $\epsilon$
is not a square.  Now assume that $W'$ occurs $n'$ times in the
decomposition of $\StGwp$ and $W''$ occurs $n''$ times.  Then, if
$\epsilon_1$ is a square and $\epsilon_2$ isn't, the above discussion
shows that
\begin{align*}
  \dim \Stgisp{\epsilon_1} - \dim \Stgisp{\epsilon_2} &= {n'}^2 +
  {n''}^2 - 2n'n'' \\
  &= (n' - n'')^2.
\end{align*}

This is a bit misleading, however, because in this case $n'$ and $n''$
are equal, so the dimension of $\Step$ is the same for all $\eps$.  We
can see this by calculating $n'$ and $n''$ using
Ligozat~\cite{ligozat}, Proposition~II.1.3.2.1: the characters of $W'$
and $W''$ only differ in matrices that are conjugate to $\smallmattwo
1*01$, and the only place that such matrices occur in the formula
given there is in the term $\sum_{a \mod p} \chi(\smallmattwo 1a01)$,
which equals $(p+1)/2$ both for $\chi = \chi_{W'}$ and $\chi =
\chi_{W''}$.

As a corollary, this implies that there are no CM-forms in $\StGwp$
for $p \equiv 1 \pmod 4$.  For if there were such a form $g$, it would
generate an irreducible representation $R_g \subset \StGwp$, all of
whose elements would be CM-forms; there would then be a form in $R_g
\otimes (R_g \circ \theta_{-1})$ that is invariant under $\PSLtp$.
But such a form would be a CM-form in $\Stgisp{-1}$, so
Theorem~\ref{thm:kernel-cm} would then imply that the dimension of
$\Step$ for $\eps$ a non-square is strictly smaller than the dimension
of $\Stgisp{-1}$, contradicting our calculations above.

Let us now turn to the case where $p \equiv 3 \pmod 4$.  This time,
$R_i \simeq \overline{R_i}$ unless $R_i \simeq X'$ or $X''$, where
$X'$ and $X''$ are the irreducible representations of dimension
$(p-1)/2$; $\overline{X'} \simeq X''$ and vice-versa.  Similarly, $R_i
\circ \theta_\epsilon \simeq R_i$ unless $R_i \simeq X'$ or $X''$ and
$\epsilon$ is not a square mod $p$; if it is, $X' \circ
\theta_\epsilon \simeq X''$ and vice-versa.  Thus, if $X'$ occurs $n'$
times and $X''$ occurs $n''$ times in the decomposition of $\StGwp$,
\begin{align*}
  \dim \Stgisp{\epsilon_1} - \dim \Stgisp{\epsilon_2} &= 2n'n'' - ({n'}^2 +
  {n''}^2)\\
  &= -(n' - n'')^2,
\end{align*}
where $\epsilon_1$ is a square mod $p$ and $\epsilon_2$ isn't.  Since
$-1$ is not a square, the dimension is maximized when $\epsilon = -1$, 
agreeing with Corollary~\ref{cor:genus-maximized}.

This time, however, $n'-n''$ is non-zero.  We can't calculate it as
easily as we calculated it in the previous case, because the method used
there calculates the number of times a representation occurs plus the
number of times that its complex conjugate occurs, and here the
character is no longer totally real.  Instead, we refer to
Hecke~\cite{hecke-xp-h}, where he proves that the difference is equal
to the class number $h(-p)$ of $\qq(\sqrt{-p})$.  Thus,
\[
\dim \Stgisp{-1}-\dim \Stgisp{1} = h(-p)^2.
\]
As before, this implies that there are exactly $h(-p)\cdot(p-1)/2$
CM-forms contained in $\StGwp$; they have been constructed by Hecke in
\cite{hecke-theta}.  We shall review his construction in
Section~\ref{sec:examples}, and use them to write down the Hecke
kernel $\Kis(p)$ explicitly.  We shall also show how to use the theory
outlined in this Section to perform explicit calculations of spaces
$\Step$ for small primes.

To recap:

\begin{theorem}
  \label{thm:kp-1}
  If $p$ is a prime congruent to 1 mod 4 then there are no CM-forms
  contained in $\StGwp$ and the Hecke kernel $\Kis(p)$ is zero.  If $p>3$
  is congruent to 3 mod 4 then there are $h(-p)\cdot (p-1)/2$ CM-forms 
  contained in $\StGwp$ and $\Kis(p)$ has dimension $(h(-p))^{2}$, where
  $h(-p)$ is the class number of $\qq(\sqrt{-p})$.
  \qed
\end{theorem}

\section{Examples}
\label{sec:examples}

\subsection*{$\Xis{-1}(7)$}

The first $\Xise(p)$ to have a non-zero $(2,2)$-cusp form is
$\Xis{-1}(7)$, as can be seen by looking at Table~1 in Kani and
Schanz~\cite{kani-modular} (and using Corollary~\ref{cor:cusp-dim-coh}
above); in fact, we see that $\dim \Stgis{-1}{7} = 1$.  We can
explicitly determine a non-zero form in this space as follows:

Conjugating $\Ga_w(7)$ by $\smallmattwo 7001$, we can think of
$X_w(7)$ as lying between the curves $X_0(49)$ and $X_1(49)$.  The
former is an elliptic curve (after choosing a base point); its L-series
gives rise to a weight two cusp form
\[
f(z) = \sum_{m>0} c_mq^m
\]
on $X_0(49)$ and $X_w(7)$.  (Here, $q = \expn z$ if we are thinking of
$f$ as a form on $X_0(49)$ and $q=\expn{z/7}$ if we are thinking of
$f$ as a form on $X_w(7)$.)  If $\chi$ is a non-trivial character on
$\zmods7$ such that $\chi(-1) = 1$ then the functions
\[
f_\chi(z) = \sum_{m>0} c_m\chi(m)q^m
\]
and
\[
f_{\chi^2}(z) = \sum_{m>0} c_m\chi^2(m)q^m
\]
are also modular forms in $S_2(\Ga_w(7))$, by Shimura~\cite{shimura},
Proposition~3.64; since the latter space is three-dimensional,
$\{f,f_\chi,f_{\chi^2}\}$ forms a basis for it.  For $n \in \zmods7$, we
have $f_\chi|_{\sigma_a} = \chi^2(a)f_\chi$ and
$f_{\chi^2}|_{\sigma_a} = \chi(a)f_{\chi^2}$.

To produce an element of $\Stgis{-1}{7}$, we have to find a form
contained in $S_2(\Ga_w(7))\otimes S_2(\Ga_w(7))$ that is fixed by
$\PSLt{7}$ (acting on the second factor via $\theta_{-1}$).  For our
form to be fixed by the matrices $\smbigparens{\sigma_a,\sigma_a}$, it
has to be of the form
\[
a_0\cdot f\otimes f+a_1\cdot f_\chi\otimes f_{\chi^2} + a_2\cdot
f_{\chi^2}\otimes f_\chi.
\]
And for our form to be fixed by the matrix $\pairsmallmattwo
1{-1}011101$, we must have $a_0 = a_1 = a_2$.  However, those constraints
leave us with only a one-dimensional space of possible cusp forms, and
since $\Stgis{-1}{7}$ is non-empty, we see that it must be generated
by the form
\[
g = \frac 13(f\otimes f + f_\chi\otimes f_{\chi^2} + f_{\chi^2}\otimes
f_\chi) = \sum_{m_1 \equiv m_2 \narpmod 7}
c_{m_1}c_{m_2}q_1^{m_1}q_2^{m_2},
\]
where the $c_i$'s are the coefficients of $f$ as above.

Now that we've got our form $g$ in hand, we'd like to relate it to
some of our general theorems about forms in $\SkeNb$.  Note that $g$
has lots of Fourier coefficients that are zero: not only is $\cmi(g)$
zero unless $m_1 \equiv m_2 \pmod 7$, but it's also zero unless the
$m_i$'s are squares mod $7$.  (This follows from the fact that the
elliptic curve $X_0(49)$ has complex multiplication by
$\qq(\sqrt{-7})$.)  By Proposition~\ref{prop:kernel-characterization},
our form is therefore in $\Kis(7)$; indeed, $\Stgis{1}{7}$ is trivial.

\subsection*{$\Xis{-1}(p)$ for $p \equiv 3 \pmod 4$}

The above may look like a general recipe for producing forms on
$\Xise(p)$ out of forms on $X_0(p^2)$, but it isn't.  To see why, note
that the transition involved two steps: matching up characters, which
involved checking invariance under the matrices
$\smbigparens{\sigma_a,\sigma_a}$, and making sure that certain
Fourier coefficients were zero, which involved checking invariance
under the matrices $\pairsmallmattwo 1\eps011101$.  Thus, we checked
that our putative form is invariant under the subgroup $B(p)$ of
upper-triangular matrices, not all of $\PSLtp$.  The reason why we
could get away with that above was that we knew a lot about
$S_2(\Ga_w(7))$ and that the dimension of $\Stgis{-1}{7}$ was 1.

Fortunately, all is not lost for more general $p$.  Let $R_1$ and
$R_2$ be an irreducible representation occurring in $\StGwp$.  As the
discussion in Section~\ref{sec:prime} showed, $R_1 \otimes R_2$
contributes 1 to the dimension of $\Step$ iff $R_1 = R_2 \circ
\theta_\epsilon$.  Now, assume that that is indeed the case, and that,
furthermore, $R_1$ is irreducible as a representation of $B(p)$.
Writing $\chi_i$ for the character of $R_i$, it will then also be the
case that
\[
\bracket{\chi_1\cdot(\chi_2\circ\theta_\eps),1_{B(p)}}_{B(p)} =
\bracket{\chi_1,\overline{\chi_2\circ\theta_\eps}}_{B(p)} = 1.
\]
But this says that there's only a one-dimensional space of vectors in
$R_1\otimes R_2$ that is fixed by $B(p)$, and since there is also a
one-dimensional space of vectors in $R_1\otimes R_2$ that is fixed by
$\PSLtp$, they must be the same space.  Thus, under the hypothesis
that our representation is irreducible when considered as a
representation of $B(p)$, we can test to see whether an element of
$R_1\otimes R_2$ is a cusp form on $\Xise(p)$ simply by making sure
that it is invariant under $\smbigparens{\sigma_n, \sigma_n}$ and
$\pairsmallmattwo 11011\eps01$.

To make this concrete, assume that $p$ is congruent to $3 \pmod 4$ but
not equal to 3 and that $\eps = -1$.  The character table for $\PSLtp$
is given in Section~\ref{sec:prime}; checking the non-trivial
characters listed there, we see that $X'$ and $X''$ remain irreducible
when restricted to $B(p)$.  Thus, if we can produce representations
isomorphic to $X'$ or $X''$ in $\StGwp$, we'll be able to explicitly
write down forms in $\Stmp$.  We saw that there should be $h(-p)$ such
representations coming from CM-forms; they would be good ones to look
for. 

\newcommand{\mubar}{{\overline{\mu}}}

Fortunately, those representations are produced in
Hecke~\cite{hecke-theta}.  They are defined as follows: let $I$ be an
integral ideal in $\qmp$ with norm $A$ and let $\rho$ be an element of
$I$.  We define a theta series as follows:
\[
\theta_H(z;\rho,I,\sqrt{-p}) = 
\sum_{\substack{\mu \in I\\\mu \equiv \rho \narpmod{I\sqrt{-p}}}}
\mu \expn{z\frac{\mu\mubar}{pA}},
\]
where $\mubar$ is the complex conjugate of $\mu$.  Letting $V_I$ be
the vector space generated by the functions
$\theta_H(z;\rho,I,\sqrt{-p})$ for $\rho \in I$, the results of
Hecke~\cite{hecke-theta}, \S 7 show that $V_I$ only depends on the
ideal class of $I$, that these $\theta_H$'s are elements of $\StGwp$,
and that $V_I$ is a representation of $\PSLtp$ isomorphic to $X'$.
This gives us our desired $h(-p)$ different copies of $X'$.

Now that we have our representations, we follow the same program as in
the $\Xis{-1}(7)$ case:

\begin{theorem}
  \label{thm:cm-forms-explicit}
  Let $p$ be a prime congruent to 3 mod 4.  For each ideal class of
  $\qq(\sqrt{-p})$, fix an integral ideal $I$ in that class and an
  element $\alpha_I$ of $I$ that's not contained in $I\sqrt{-p}$.  Let
  \[
  f_I = \sum_{a\in \zmods p} \theta_H(z;a\jacobi ap\alpha_I,I,\sqrt{-p})
  \]
  have the Fourier expansion
  \[
  f_I(z) = \sum_{m>0} c_{I,m} q^m,
  \]
  where $q = \expn{z/p}$.  If $I_1$ and $I_2$ are (not necessarily
  distinct) ideal classes then the function
  \[
  f_{I_1,I_2}(z_1,z_2) = \sum_{m_1 \equiv m_2 \narpmod p}
  c_{I_1,m_1}c_{I_2,m_2}q_1^{m_1}q_2^{m_2}.
  \]
  is an element of $\Stmp$ contained in $\Kis(p)$; furthermore, the
  $f_{I_1,I_2}$'s give a basis for $\Kis(p)$ as $I_1$ and $I_2$ vary over 
  the ideal classes of $\qq(\sqrt{-p})$.
\end{theorem}

\begin{proof}
  The same argument as in the $p=7$ case shows that multiples of
  $f_{I_1,I_2}$ are the only elements of $V_{I_1}\otimes V_{I_2}$
  invariant under $B(p)$, so they are indeed elements of $\Step$.
  Assuming that we can show that they are in $\Kis(p)$,
  Theorem~\ref{thm:kp-1} shows that they give us a basis.  Thus, by
  Theorem~\ref{thm:kernel-cm}, we just have to verify that the forms
  $f_I$ are CM-forms.
  
  This can be seen as follows: by definition,
  \[
  c_m(\theta_H(z;\rho,I,\sqrt{-p}))
  = \sum_{\substack{\mu \in I\\\mu\equiv \rho \narpmod{I\sqrt{-p}} \\
      \mu\mubar = mA}} \mu,
  \]
  where $A$ is the norm of $I$ But $\mu\mubar$ is a square mod $p$ for
  all $\mu$ in the ring of integers of $\qq(\sqrt{-p})$, as is $A$, so
  $c_m$ is zero unless $m$ is a square mod $p$.  Thus, $f_I$ is
  invariant under twisting by the quadratic character of $\zmods p$,
  hence a CM-form.
\end{proof}

\end{document}